\documentclass{elsart}

\usepackage{amsmath,amsfonts,amssymb,epsfig}
\usepackage{graphics,graphicx}
\usepackage{color}
\usepackage{epsf,epsfig,psfrag}
\def\epsfx#1#2{\leavevmode\epsfxsize=#1 \epsfbox{#2}}

\usepackage{t1enc}

\usepackage[english]{babel}

\journal{Stochastic Processes and their Applications}

\def\epsfx#1#2{\leavevmode\epsfxsize=#1 \epsfbox{#2}}

\def \cC{{\mathcal C}}

\def \cD{{\mathcal D}}

\def \N{{\mathbf N}}

\def \cal{\mathcal}

\def\tha{{\tau_{\Delta^\alpha}}}

\def\thb{\tau_{\Delta^\beta}}

\def\tro{\tau^{r_0}}

\def\H{{\mathbb{H}}}

\def\cA{{\mathcal{A}}}

\def\cD{{\mathcal{D}}}

\def\dcD{ \partial{\mathcal{D}} }

\newcommand{\pd}{\mathbf{pd}}

\def\qed{
$\ \Box$
}
\renewcommand{\theequation}{\thesection.\arabic{equation}}
\renewcommand{\thesection}{\arabic{section}}
\newcommand{\mysection}{\setcounter{equation}{0} \section}

\newcommand{\1}{\mathbf{1}}

\newcommand{\I}{\1}

\newcommand{\E}{\mathbb{E}}

\newcommand{\p}{\mathbb{P}}

\renewcommand{\P}{\p}

\newcommand{\R}{\mathbb{R}}

\newcommand{\D}{\mathcal{D}}

\renewcommand{\H}{{\mathbf H}}

\newcommand{\deD}{\partial D}

\newcommand{\PD}{\mathcal{PD}}

\def \pdti{{\rm \pi}_{\bar D_{t_i}}}

\def \pdtip{{\rm \pi}_{\bar D_{t_{i+1}}}}

\def \pddt{{\rm \pi}_{\deD_t}}

\def \={\stackrel{\E}{=}}

\def \<{\stackrel{\E}{\le}}

\newcommand{\No}{\mathcal {N}}

\newcommand{\F}{\mathcal{F}}

\newcommand \A[1]{{\bf (#1)}}

\def \etal{{\it et al. }}

\def \Tr {{\rm Tr }}

\newcommand \bfrac[2]{{\displaystyle \frac{#1}{#2}}}

\def\bint#1^#2{\displaystyle{\int_{#1}^{#2}}}

\def\bsum#1^#2{\displaystyle{\sum_{#1}^{#2}}}

\newcommand \Err[1]{{\rm Err}_{#1}(T,\Delta,g,f,k,x)}

\newcommand \ellErr[1]{{\rm Err}_{#1}(\Delta,g,f,k,x)}

\def \tautx {\tau^{t,x}}

\def \xtx{X^{t,x}}

\begin{document}

\begin{frontmatter}



\title{Stopped diffusion processes: boundary corrections and overshoot}


\author[Grenoble]{Emmanuel Gobet\corauthref{cor1}},
  \ead{emmanuel.gobet@imag.fr}
  \author[P7]{St\'ephane Menozzi}
  \ead{menozzi@math.jussieu.fr}
  \address[Grenoble]{Laboratoire Jean Kuntzmann, Universit{\'e} de Grenoble and CNRS, BP 53, 38041 Grenoble Cedex 9, FRANCE}
\address[P7]{Laboratoire de Probabilit\'es et Mod\`eles Al\'eatoires, Universit\'e Denis Diderot Paris 7, 175 rue de Chevaleret 75013 Paris, FRANCE.}
\corauth[cor1]{Corresponding author}
\begin{abstract}
For a stopped diffusion process in a multidimensional time-dependent domain $\D$, we propose and analyse a new procedure consisting in simulating the process with an Euler scheme with step size $\Delta$ and stopping it at discrete times $(i\Delta)_{i\in\N^*}$ in a modified domain, whose boundary has been appropriately shifted. The shift is locally in the direction of the inward normal $n(t,x)$ at any point $(t,x)$ on the parabolic boundary of $\D$, and its amplitude is equal to $0.5826 (...) |n^*\sigma|(t,x)\sqrt \Delta$ where $\sigma$ stands for the diffusion coefficient of the process. The procedure is thus extremely easy to use. In addition, we prove that the rate of convergence w.r.t. $\Delta$ for the associated weak error is higher than without shifting, generalizing previous results by \cite{broa:glas:kou:97} obtained for the one dimensional Brownian motion. For this, we establish in full generality the asymptotics of the triplet exit time/exit position/overshoot for the discretely stopped Euler scheme. Here, the overshoot means the distance to the boundary of the process when it exits the domain. Numerical experiments support these results.
\end{abstract}

\end{frontmatter}

\section{Introduction}

\subsection{Statement of the problem}

\begin{figure}[h]  
\centerline{\epsfx{6cm}{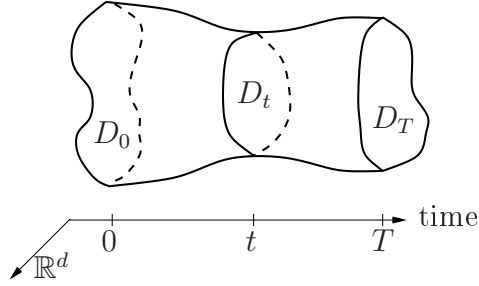}}  
\caption{Time space domain and its time-sections.\label{fig1bis}} 
\end{figure}

We consider a $d$-dimensional diffusion process whose dynamics is given by
\begin{equation}
\label{diff}
X_t=x+\bint{0}^{t}b(s,X_s)ds+\bint{0}^{t}\sigma(s,X_s)dW_s
\end{equation}
where  $W$ is a standard $d'$-dimensional Brownian motion defined on a
filtered probability space $(\Omega,\F,(\F_t)_{t\ge 0},\P) $
satisfying the usual conditions. The mappings $b$ and $\sigma $ are
Lipschitz continuous in space and locally bounded in time, so that
\eqref{diff} has a unique strong solution. 
We consider $(D_t)_{t\ge 0}$, a time-dependent family of smooth bounded domains of $\R^d$, that is also smooth with respect to $t$ (we refer to paragraph \ref{time domains} for a precise definition). See Figure \ref{fig1bis}. For a fixed deterministic time $T>0$, this defines a time-space domain $$\D=\bigcup_{0<t<T}\{t\}\times D_t=\{(t,x): 0<t<T, x\in D_t\}\subset]0,T[\times\R^d.$$
Cylindrical domains are specific cases of time-dependent domains of the form $\D=]0,T[\times D$, where $D$ is a usual domain of $\R^d$ ($D_t=D$ for any $t$). Time-dependent domains in dimension $d=1$ are typically of the form $\D=\{(t,x):0<t<T, \varphi_1(t)<x<\varphi_2(t)\}$ for two functions $\varphi_1$ and $\varphi_2$ (the time-varying boundaries).

\psfrag{1bD0}{ $D_0$}  
\psfrag{1bDt}{$D_t$}  
\psfrag{1bDT}{$D_T$}  
\psfrag{1btime}{time}  
\psfrag{1bt}{$t$}  
\psfrag{1b0}{$0$}  
\psfrag{1bT}{$T$} 
\psfrag{1bRd}{$\R^d$} 
Now, set $\tau:=\inf\{t> 0:X_t\not \in D_t \}$, then $\tau\land T$ is the first exit time of $(s,X_s)_s$ from the time-space domain $\D$.
Given continuous functions $g,f,k:  \bar \cD \rightarrow \R$, we are
interested in estimating the quantity
\begin{align}
\label{quant}
\hspace{-3mm}\E_x[g(\tau\land T,X_{ \tau\land T})Z_{\tau\land T}+\int_{0}^{\tau\land T}Z_sf(s,X_s)ds],
Z_s=\exp(-\bint{0}^{s}k(r,X_r)dr),
\end{align}
where as usual $\E_x[.]:=\E[.|X_0=x]$ (resp. $\P_x[.]:=\P[.|X_0=x]$).
The approximation of such quantities is a well known issue in
finance, since it represents in this framework the price of a
barrier option, see e.g. Andersen and Brotherton-Ratcliffe
\cite{ande:brot:96}. These quantities also arise through the
Feynman-Kac representation of the solution of a parabolic PDE with
Cauchy-Dirichlet boundary conditions, see Costantini \etal
\cite{cost:elka:gobe:06}. They can therefore also be related to
problems of heat diffusion in time-dependent domains. 

We then choose to approximate the expectation in \eqref{quant} by
Monte Carlo simulation. This approach is natural and especially
relevant compared to deterministic methods if the dimension $d$ is large. To this end we
approximate the diffusion \eqref{diff} by its Euler scheme with time
step $\Delta>0$ and discretization times $(t_i=i\Delta=iT/m)_{i\ge 0}$ ($m\in\N^*$ so that $t_m=T$). For $t\ge 0$, define $\phi(t)=t_i$ for $t_i\le t<t_{i+1}$ and
introduce
\begin{align}
\label{euler}
X_t^\Delta=x+\bint{0}^{t}b(\phi(s),X_{\phi(s)}^\Delta)ds+\bint{0}^{t}\sigma(\phi(s),X_{\phi(s)}^\Delta )dW_s.
\end{align}
We now associate to \eqref{euler} 
the discrete exit time
$\tau^\Delta:=\inf\{t_i>0:X_{t_i}^\Delta\notin
D_{t_i}\}$.
Approximating the functional
$
V_\tau:=g(\tau\land T,X_{\tau\land T}
)Z_{\tau\land T}+\int_0^{\tau\land T}Z_sf(s,X_s)$ $ds$
by
\begin{align*}
& 
V_{\tau^\Delta}^\Delta
:= g(\tau^{\Delta}\land T,
X^{\Delta}_{\tau^{\Delta}\land T}
)Z^{\Delta}_{\tau^{\Delta}\land T}+\int_{0}^{\tau^{\Delta}\land T} Z^{\Delta}_{\phi(s)} f(\phi(s),X^{\Delta}_{\phi(s)}) ds\\
&\quad\text{ with} \quad Z^{\Delta}_{t}=e^{-\int_0^{
t
} k(\phi(r),X^{\Delta}_{\phi(r)}) dr},
\end{align*}
we introduce the quantity
\begin{align}
\label{erreur}
\Err{} =\E_x[V_{\tau^\Delta}^\Delta-V_\tau]
\end{align}
that will be referred to as the weak error. 

Note that in $V_{\tau^\Delta}^\Delta $, on $\{\tau^\Delta\le T\}$ $g$ is a.s. not evaluated on the side part $\bigcup_{0\le t\le T}\{t\}\times\partial D_t$ of the boundary ($g$ must be understood as a function defined in a neighborhood of the boundary). At first sight, this approximation can seem coarse. Anyhow, it does not affect the convergence rate and really reduces the computational cost with respect to the alternative that would consist in taking the projection on $\dcD$. It is a commonly observed phenomenon that the error is positive when $g$ is positive (overestimation of $\E_x(V_\tau)$), because we neglect the possible exits between two discrete times: see Boyle and Lau \cite{boyl:lau:94}, Baldi \cite{bald:95}, Gobet and Menozzi \cite{gobe:meno:spa:04}. In addition, it is known that the error is of order $\Delta^{1/2}$: see \cite{gobe:meno:spa:04} for lower bound results, see \cite{gobe:meno:semi:07} for upper bounds in the more general case of Itô processes. But so far, the derivation of an error expansion $\E_x[V_{\tau^\Delta}^\Delta-V_\tau]=C\sqrt \Delta +o(\sqrt \Delta)$ had not been established: this is one of the intermediary results of the current work (see Theorem \ref{dev}).

\psfrag{1partialDt}{\hspace{-4mm}$\partial D_t$}  
\psfrag{1Ddeltat}{\ $D^\Delta_t$}  
\begin{figure}[h]  
\centerline{\epsfx{4cm}{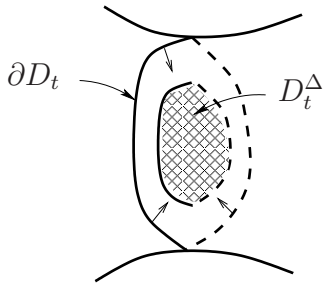}}  
\caption{The boundary $\partial D_t$ and the smaller domain $D^\Delta_t$.\label{fig1}}
\end{figure}  
Our goal goes beyond this result, by designing a simple and very efficient improved procedure. We propose to stop the Euler scheme at its exit of a smaller domain in order to compensate the underestimation of exits and to achieve an error of order $o(\sqrt \Delta)$. The smaller domain is defined by its time-section
$$D^\Delta_t=\{x\in D_t: {\bf d}(x,\deD_t)> c_0 \sqrt \Delta |n^*\sigma(t,x)|\}$$
where $n(t,x)$ is the inward normal vector at the closest point of $x$ on the boundary $\deD_t$, see Figures \ref{fig1} and \ref{fig2} for details\footnote{the closest point of $x$ may not be unique for points $x$ far from $\deD_t$. But since the above definition of $D^\Delta_t$ involves only points close to the boundary, this does not make any difference.}. We shall interpret $|n^*\sigma(t,x)|$ as the noise amplitude along the normal direction to the boundary. The constant $c_0$ is defined later in (\ref{eq:c0}) and equals approximatively $0.5826(\ldots)$. Thus, the associated exit time of the Euler scheme is given by 
$$\hat\tau^{\Delta}=\inf\{t_i>0: X^\Delta_{t_i}\not \in D_{t_i}^\Delta\}\le \tau^\Delta.$$ 
The new Monte Carlo scheme consists in simulating independent realizations of 
\begin{align*}
 V_{\hat \tau^\Delta}^\Delta=g(\hat\tau^{\Delta}\land T,X^{\Delta}_{\hat\tau^{\Delta}\land T})
Z^{\Delta}_{\hat\tau^{\Delta}\land T}+\int_{0}^{\hat\tau^{\Delta}\land T} Z^{\Delta}_{\phi(s)} f(\phi(s),X^{\Delta}_{\phi(s)}) ds
\end{align*}
and averaging them out to get an estimator of the required quantity $\E_x(V_\tau)$. Our main result (Theorem \ref{corr}) is that the asymptotic bias w.r.t. $\Delta$ is significantly improved:
$$\E_x[V_{\hat \tau^\Delta}^\Delta-V_\tau]=o(\sqrt \Delta)$$
(instead of $C\sqrt \Delta +o(\sqrt \Delta)$ before). This improvement has been already established in the case of the one-dimensional Brownian motion \cite{broa:glas:kou:97} in the context of computational finance, exploiting heavily the connection with Gaussian random walks and some explicit computations available in the Brownian motion case. 

\subsection{Contribution of the paper}
To achieve the results in the current very general framework, we combine several ingredients (which correspond to the main steps of the proofs).
\begin{enumerate}
\item We first expand the error $\E_x[V_{\tau^\Delta}^\Delta-V_\tau]$ related to the use of the discrete Euler scheme in the domain $\D$. Although this issue deserved many studies in the literature, the expansion results are new. We prove that it relies on the study of the weak convergence of the triplet (exit time, position at exit time, renormalized overshoot at exit time), that is $(\tau^\Delta, X^\Delta_{\tau^\Delta}, \break \Delta^{-1/2}d(X^\Delta_{\tau^\Delta}, \deD_{\tau^\Delta}))$, as $\Delta$ goes to 0. This weak convergence result is crucial in this work and it is new (see Theorem \ref{loi_limite}).\\ Then, combining this with sharp techniques of error analysis, we derive an expansion of the form
$\Err{}=C\sqrt \Delta+o(\Delta)$ in the very general framework of
stopped diffusions in time-dependent domains.
\item Second, we analyse the impact of the boundary shifting, in the continuous time problem (see paragraph~\ref{par:boundary:correction}). This is related to the differentiability of $\E_x(V_\tau)$ w.r.t. the boundary and it has been addressed in \cite{cost:elka:gobe:06}. We apply directly their results. Then, we obtain the global error estimate of the boundary correction procedure (Theorem \ref{corr}).
\end{enumerate}



We mention that the
previous results about the error expansion and correction still hold
in the stationary setting, see Section \ref{stat}, which also seems to be new. A numerical application is discussed in Section \ref{resnum}. Complementary tests are presented in \cite{gobe:09}, showing that the boundary correction procedure is very generic and seems to work without Markovian property for $X$. This feature will be investigated in further research.

Let us finally mention that we could also consider the diffusion process discretely stopped: expansion and correction results below would remain the same.

\subsection{Comparison with results in literature}
Up to now, the behavior of \eqref{erreur} had mainly been analysed
for cylindrical domains, in the killed case, without source and
potential terms (i.e. when the error
writes
${\rm Err}(T,\Delta,g,0,0,x)=\E_x[g(X_T^\Delta)\I_{\tau^\Delta>T}]-\E_x[g(X_T)\I_{\tau>T}]$). Let
us first mention the work of Broadie \etal \cite{broa:glas:kou:97}, who first derived the boundary shifting procedure in the one
dimensional geometric
Brownian motion setting (Black and Scholes model). In \cite{gobe:00}
and \cite{gobe:meno:spa:04}, it had been shown that, under some (hypo)ellipticity conditions on the
coefficients and some smoothness of the domain and the coefficients, ${\rm Err}(T,\Delta,g,0,0,x)$
was lower and upper bounded at order $1/2$ w.r.t. the time-step
$\Delta$. Also, an expansion result for the killed Brownian motion in
a cone as well as the associated correction procedure are available in
\cite{meno:06}.

All these works emphasize that the crucial quantity to
analyse in order to obtain an expansion is the overshoot above
the spatial boundary of the discrete process. 
In the Brownian one-dimensional framework such analysis goes back to
Siegmund \cite{sieg:79} and Siegmund and Yuh \cite{sieg:yuh:82}. Also a non linear renewal theory for random walk, i.e. for a curved boundary, had been developed by Siegmund and \textit{al.}, see \cite{sieg:85} and references therein, Woodroofe \cite{wood:82} and Zhang \cite{zhan:88}. We
manage to extend their results to obtain the asymptotic distribution
of the overshoot of the Euler scheme, see Sections \ref{main_results}
and \ref{proof_res_os}. Concerning the asymptotics of the overshoot of
stochastic processes, let us mention the works of Alsmeyer \cite{alsm:94} or Fuh and Lai \cite{fuh:lai:01} for ergodic Markov chains and Doney and Kyprianou for Lévy processes \cite{done:kypr:06}. These works are all based on renewal arguments.



Finally, for simulating stopped diffusions we also mention the alternative technique based on Random Walks on Spheres. This method allows to derive 
a bound for the weak error associated to the approximation of $\E[V_\tau]$ in the elliptic setting for a cylindrical domain, see Milstein \cite{mils:97}. The same approach has also been exploited to obtain some strong error or pathwise bounds for a bounded time-space cylindrical domain, see Milstein and Tretyakov \cite{mils:tret:99}. Recently, Deaconu and Lejay \cite{deac:leja:06} have developed similar algorithms, but based on random walks on rectangles.
However, computationally speaking, our approach is presumably more direct.

\subsection{Outline of the paper}
Notations and assumptions used throughout the paper are stated in
Section \ref{not_and_ass}. In Section \ref{main_results} we give 
our main results concerning the asymptotics of the overshoot, the error expansion and the boundary
correction. 
These results are proved in Section
\ref{proof_res_os}, which is the technical core of the
paper. Eventually, Section \ref{stat} deals with the stationary extension of our
results. We still manage to obtain an expansion and a correction for
elliptic PDEs. Some technical results are postponed to the Appendix.

\subsection{General notation and assumptions}
\label{not_and_ass}

\subsubsection{Miscellaneous}\hspace*{10pt}\\
\noindent $\bullet$ {\em Differentiation.}
For smooth functions $g(t,x)$, we denote by $\partial^{\beta}_xg(
t,x)$ the derivative
of $g$ w.r.t. $x$ according to the
multi-index $\beta$, whereas the time-derivative of $g$ is denoted
by $\partial_tg(t,x)$. The notation $\nabla
g(t,x)$ stands for the usual
gradient w.r.t. $x$ (as a row vector) and the Hessian matrix of $
g$ (w.r.t. the space variable $x$)
is denoted by $Hg(t,x)$.\\
The second order linear operator $L_t$ below stands for the infinitesimal
generator of the diffusion process $X$ in \eqref{diff} at time $t$ :
\begin{equation}\label{eq:generator}L_tg(t,x)=\nabla g(t,
x)b(t,x)+\frac 12\Tr(Hg(t,x)[\sigma\sigma^{*}](t,x)).\end{equation}

\noindent $\bullet$ {\em Metric.} 
The Euclidean norm is denoted by $|\cdot|$.\\ 
We set
$B_d(x,\epsilon )$ for the usual Euclidean $d$-dimensional open ball with
center $x$ and radius $\epsilon$ and ${\bf d}(x,C)$ for the Euclidean
distance of a point $x$ to a closed set $C$. 
The $r$-neighborhood of $C$ is denoted by $V_{C}(r)=\{x: {\bf d}(x,C)\le r \}$ ($r\ge 0$).
\\
\noindent$\bullet$ {\em Functions.\/} For an open set $\D'\subset
\R\times\R^d$ and $l\in\mathbf N$, $\mathcal C^{\lfloor\frac l2
\rfloor ,l}(\D')$ (resp. $\mathcal C^{\lfloor\frac l2\rfloor
,l}(\overline {\D'})$) is the space of continuous functions $ f$
defined on $\D'$ with continuous derivatives
$\partial^{\beta}_x\partial^j_tf$ for $|\beta |+2j\leq l$ (resp.
defined in a neighborhood of $\overline {\D'}$). 
Also, for $a=l+\theta,\ \theta\in ]0,1],\ l\in\N$, we denote by
$\H_a(\cD')$ (resp. $\H_a(\bar{\cD'}) $) the Banach space of
functions of $\cC^{\lfloor \frac l 2\rfloor ,l}(\cD') $ (resp.
$\cC^{\lfloor \frac l 2\rfloor ,l}(\bar \cD')$ ) having $l^{\rm th}$
space derivatives uniformly $\theta$-Hölder continuous and
$\lfloor l/2\rfloor $ time-derivatives uniformly $ (a/2-\lfloor  l /2\rfloor)
$-Hölder continuous, see Lieberman \cite{lieb:96}, p. 46 for
details. We may simply write $\cC^{\lfloor \frac l 2\rfloor ,l}$ or
${\H}_{a}$ when  $\D'=\R\times\R^d$.
\\
\noindent$\bullet$ {\em Floating constants}. As usual, we use the same
symbol $C$ for all finite, non-negative constants which appear in our computations : they may depend on
$\cD,T, b,\sigma, g,f,k $ but they will not depend on $\Delta$ or
$x$. We reserve the notation $c$ for constants also independent of
$T$, $ g,f$ and $k$. Other possible dependences will
be explicitly indicated. \\
In the following
$O_{pol}(\Delta)$ (resp. $O(\Delta)$) stands for every quantity
$R(\Delta)$ such that, for any  $k\in \N$ one has $|R(\Delta)|\leq C_k\Delta^k$ (resp. $|R(\Delta)|\leq C\Delta$) for a constant $C_k>0$ 
(uniformly in the starting point $x$).

\subsubsection{Time-space domains}\label{time domains}
Below, we introduce some usual notations for such domains (see e.g. \cite{frie:64}, \cite{lieb:96}). In what follows, for any $t\ge 0$, $D_t$ is a non empty bounded domain of $\R^d$, that coincides with the interior of its closure (see \cite{frie:64}, Section 3.2). We then define the time-space domain by $\D:=\bigcup_{0<t<T}\{t\}\times D_t\subset]0,T[\times\R^d$, see Figure \ref{fig1bis}.

Regularity assumptions on the domain $\D$ will be formulated in
terms of H\"older spaces with {\em time-space\/}  variables (see
\cite{lieb:96} p.46 and \cite{frie:64} Section 3.2). Namely, we say
that the domain $\cD$ is of class $\H_a,\ a\ge 1$ if for every
boundary point $(t_0,x_0)\in \bigcup_{0\le t \le T}\{t\}\times \deD_t$,
there exists a neighborhood $]t_0-\varepsilon_0^2
,
t_0+\varepsilon_0^2[\times B_d(x_0,\varepsilon_0)$, an index $1\le i\le d$ and a function $\varphi_0\in \H_a(]t_0-\varepsilon_0^2
,t_0+\varepsilon_0^2[\times
B_{d-1}((x_0^1,...,x_{0}^{i-1},x_0^{i+1},...,x_0^d),\varepsilon_0) $
s.t.
\begin{align*}
&\big\{\cup_{0\le t \le T}\{t\}\times \deD_t\big\} \cap \big\{]t_0-\varepsilon_0^2
, t_0+\varepsilon_0^2[\times B_d(x_0,\varepsilon_0)\big\}\\
& \hspace*{15pt}:=\{(t,x)\in (]t_0-\varepsilon_0^2
,t_0+\varepsilon_0^2 [\cap[0,T])\times B_d(x_0,\varepsilon_0):\\
& \hspace*{55pt} x_i=\varphi_0(t,x_1,...,x_{i-1},x_{i+1},...,x_d)\}.
\end{align*}
If $\D$ is of class $\H{_2}$, all domains $ D_t$, for $t\in [0,T]$, satisfy the
uniform interior and exterior sphere condition with the same radius $r_0>0$. Moreover, the signed spatial distance $F$, given by
\[F(t,x)=\left\{\begin{array}{ll}
-{\bf d}(x,\partial D_t),&\ \mbox{\rm for }x\in D_t^c,\,{\bf d}
(x,\partial D_t)\le r_0,\,0\le t\le T,\\
{\bf d}(x,\partial D_t),&\ \mbox{\rm for }x\in D_t,\,{\bf d}(x
,\partial D_t)\le r_0,\,0\le t \le T,\end{array} \right.\] belongs to
$\H{_2}\left(\{(t,x):0\le t \le T,\,{\bf d}(x,\partial D_t)<r_ 0\}\right)$ (see \cite{lieb:96}, Section X.3)
and $n(t,x)=[\nabla F]^*(t,x)$ is the unit inward normal vector to $D_t $ at $\pddt(x)$
the nearest point to $x$ in $\partial D_t$ (see Figure \ref{fig2}). The function $F$ can be extended as a $\H_2([0,T]\times\R^d) $ function, preserving the sign (see \cite{lieb:96}, Section X.3).


\subsubsection{Diffusion processes stopped at the boundary}\label{section:dirichlet}

We specify the properties of the coefficients $(b,\sigma)$ in
\eqref{diff} with assumption
\begin{itemize}
\item[\A{A$_{\theta}$}] (with $\theta\in ]0,1]$)
  \begin{enumerate}
  \item[\bf 1.] {\em Smoothness}. The functions $b$ and $\sigma$ are in $\H_{1+\theta}$.
  \item[\bf 2.] {\em Uniform ellipticity}. For some $a_0>0$, it holds
$\xi^{*}[\sigma\sigma^{*}](t,x)\xi\geq a_0|\xi |^2$ for any $(t,x
,\xi )\in [0,T]\times\R^d\times\R^d$.
\end{enumerate}
\end{itemize}
We mention that the additional smoothness of $b$ and $\sigma$ w.r.t.
the time variable is required for the connection with PDEs. We also introduce assumption \A{${\mathbf A}_\theta^{'}$} for which  \textbf{2.}
is replaced by the weaker assumption
\begin{enumerate}
\item[\bf 2'.]{\em Uniform non characteristic boundary}. For some
  $r_0>0$ there exists  $a_0>0$ s.t. $\nabla F(t,x)[\sigma\sigma^{*}](t,x)\nabla F(t,x)^*\geq a_0 $
for any $(t,x)\in \bigcup_{0\le t \le T}\{t\}\times V_{\deD_t}(r_0)$.
\end{enumerate}
The asymptotic results concerning the overshoot hold true under
\A{${\mathbf A}_\theta^{'}$}, see Section \ref{mr_os}. In the following we use the superscript $t,x$ to indicate the usual Markovian dependence, i.e. $\forall s\ge t,\ X_s^{t,x}=x+\int_{t}^{s}b(u,\xtx_u)du+\int_{t}^{s}\sigma(u,\xtx_u)dW_u $.
Now let
\begin{equation}\label{eq:tps:sortie}\tautx:=\inf\{s>t:\xtx_s
\notin D_s\}\end{equation}
be the first exit time of $\xtx_s$ from $ D_s$. For functionals of the process $X$
stopped at the exit from $\D$, of the form
\begin{align}u(t,x)=&\E\big[g(\tautx \land T,\xtx_{\tautx\land T}
)e^{-\int_t^{\tautx\land T}k(r,\xtx_r)dr}\nonumber\\
&\qquad+\int_t^{\tautx\land T}e^{-\int_t^sk(r,\xtx_r
)dr}f(s,\xtx_s)ds\big],\label{eq:def:u}\end{align}
we now recall (see \cite{cost:elka:gobe:06}) that the Feynman-Kac
representation holds in the time-space domain. Introduce the parabolic
boundary $\PD=\partial\D\backslash[\{0\}\times D_0]$.
\begin{prop}{\bf [Feynman-Kac's formula and a priori estimates on} $
  u${\bf ]} \label{prop:feynman:kac:stopped}\\
  Assume \A{A$_{\theta}$}, $\D\in\H{_1}$, $k\in\H{_{\theta}}$, $f
\in\H{_{\theta}}$ and $g\in {\cal C}^{0,0}$ with $\theta\in ]0,1[$.
  Then, there is a unique solution in ${\cal C}^{1,2}(\D)\cap
{\cal C}^{0,0}(\overline{{\D}})$ to
\begin{equation}\label{eq:edp}\left\{\begin{array}{cc}
\partial_tu+L_tu-ku+f=0&{\rm i}{\rm n}\ \D,\\
u=g&{\rm o}{\rm n}\ \PD,\end{array}
\right.\end{equation}
and it is given by \eqref{eq:def:u}.\\
In addition, if for some $\theta\in ]0,1[$, $\D$ is of class $\H_{1+\theta}$, $g\in\H_{1+\theta} $ then
$u\in\H_{1+\theta}$. In particular $\nabla u$ exists and is $\theta $-Hölder continuous up to the boundary.\\
Eventually, for $\cD\in \H{_{3+\theta}}$, $k,f\in \H_{1+\theta}$,
$g\in\H{_{3+\theta}}$ satisfying the first order compatibility
condition $(\partial_t+L_T-k)g(T,x)+f(T,x)|_{x\in\dcD_T}=0 $, then the function $
u$ belongs to
$\H{_{3+\theta}}$.
\end{prop}
\textit{Proof.} The first two existence and uniqueness result for \eqref{eq:edp} are respectively implied by Theorems
5.9 and 5.10 and Theorem 6.45 in Lieberman, \cite{lieb:96}. The probabilistic
representation is then a usual verification argument, see
e.g. Appendix B.1 in \cite{cost:elka:gobe:06}. The additional
smoothness can be derived from exercise 4.5 Chapter IV in \cite{lieb:96} or Theorem 12, Chapter 3 in \cite{frie:64}.\qed

\mysection{Main Results}
\label{main_results}
\subsection{Controls concerning the overshoot}
\label{mr_os}

The overshoot is the distance of the discretely killed process to the boundary, when it exits the domain by its side. To be precise, we use $F$ the signed distance function and we consider the quantity $F(t_i,X^\Delta_{t_i})$. It remains positive for $t_i<\tau^\Delta$, and at time $t_i=\tau^\Delta$, it becomes non positive. Additionally, under the ellipticity assumption, the above inequality is strict: $F(\tau^\Delta,X^\Delta_{\tau^\Delta})<0$ a.s..
The overshoot is thus defined by $F^-(\tau^\Delta,X_{\tau^\Delta}^\Delta )$. 
Also, since $F$ is in $\H_2$ (and therefore Lipschitz continuous in time and space), it is easy to see that $F^-(\tau^\Delta,X_{\tau^\Delta}^\Delta )$ is of order $\sqrt \Delta$ (in $L_p$-norm for instance). Thus, it is natural to study the asymptotics of the rescaled overshoot
$$\Delta^{-1/2} F^-(\tau^\Delta,X_{\tau^\Delta}^\Delta ).$$

Adapting the proof of Proposition 6 in \cite{gobe:meno:spa:04} to our
time-dependent context, see also the proof of Proposition \ref{tension_STAT} for a simpler version, one has the following proposition.
\begin{prop}[Tightness of the overshoot]
\label{tension}
Assume \A{${\mathbf A}_\theta^{'} $}, and that $\cD$ is of class $\H_2 $. Then, for some $c>0$ one has
$$\sup_{\Delta >0,s\in[0,T]}\E_x[\exp(c[\Delta^{-1/2} F^-(s\wedge \tau^\Delta,X_{s\wedge\tau^\Delta}^\Delta )]^2)]<+\infty.$$
\end{prop}
It is quite plain to prove, by pathwise convergence of $X^\Delta$ towards $X$ on compact sets, that $(\tau^\Delta\land T,X_{\tau^\Delta\land T}^\Delta) $ converges in probability to $(\tau\land T,X_{\tau\land T}) $. The next theorem also includes the rescaled overshoot.
\begin{thm}[Joint limit laws associated to the overshoot] \label{loi_limite}
Assume \A{${\mathbf A}_\theta^{'} $}, and that $\cD$ is of class
$\H_2 $. Let $\varphi$ be a continuous function with compact
support. For all $t\in[0,T],\ x\in D_0,\ y\ge 0$,
\begin{align*}
&&\E_x[\I_{\tau^\Delta\le
t}Z_{\tau^\Delta}^\Delta\varphi(X_{\tau^\Delta}^\Delta)\I_{F^-(\tau^\Delta,X_{\tau^\Delta}^\Delta)\ge
y\sqrt
  \Delta}]\underset{\Delta \rightarrow 0}{\longrightarrow}\\
 &&\E_x\bigl[\I_{\tau\le t}Z_\tau\varphi(X_\tau)\bigl(1-H(y/
|  \nabla F \sigma
  (\tau,X_\tau)
|)\bigr)\bigr]
\end{align*}
with $H(y):=(\E_0[s_{\tau^+}])^{-1}\int_{0}^{y}\P_0[s_{\tau^+}>z] dz$
  and $s_0:=0, \forall n\geq 1,  s_n:=\sum_{i=1}^{n}G^i$, the $G^i$ being
i.i.d.$\ $ standard centered normal variables, $\tau^+:=\inf\{n\ge  0:
  s_n>0  \} $.
\end{thm}
In other words, $(\tau^\Delta, X_{\tau^\Delta}^\Delta, \Delta^{-1/2} F^-(\tau^\Delta,X_{\tau^\Delta}^\Delta))$ weakly converges to \break $(\tau, X_{\tau}, | \nabla F \sigma
  (\tau,X_\tau) | Y)$ where $Y$ is a random variable independent of $(\tau,X_\tau)$, and which cumulative function is equal to $H$. Actually, $Y$ has the asymptotic law of the renormalized Brownian overshoot. In the following analysis, the mean of the overshoot is an important quantity and it is worth noting that one has $\E(Y)=\frac{\E_0[s_{\tau^+}^2]}{2\E_0[s_{\tau^+}]}:=c_0$. One knows from \cite{sieg:79} that
  \begin{equation}
    \label{eq:c0}
    c_0=-\frac{\zeta(1/2)}{\sqrt{2\pi}}=0.5826...
  \end{equation}
The above theorem is the crucial tool in the derivation of our
main results. The proof is given in Section \ref{ref_pre_ll}.

\subsection{Error expansion and boundary correction}
For notational convenience introduce for $x\in D_0 $,
\begin{align*}
u(\D)&=\E_x(g(\tau\land T,X_{\tau\land T})Z_{\tau\land T}+\int_0^{\tau\land T}Z_sf(s,X_s)ds),\\
u^\Delta(\D)&=\E_x(g(\tau^{\Delta}\land T,
X^{\Delta}_{\tau^{\Delta}\land T})Z^{\Delta}_{\tau^{\Delta}\land T}+\int_{0}^{\tau^{\Delta}\land T}
  Z^{\Delta}_{\phi(s)} f(\phi(s),X^{\Delta}_{\phi(s)}) ds).
\end{align*}
\begin{thm}[First order expansion]
\label{dev}
Under \A{A$_\theta$}, for a domain of class $\H_2$, $g\in\H_{1+\theta} $,
$k,f \in\H_{1+\theta}$
and for $\Delta$
small enough
\begin{align*}&
\Err{}=u^\Delta(\D)-u(\D)\\
&=c_0 \sqrt \Delta \E_x(\1_{\tau\le T}Z_\tau (\nabla
u-\nabla g )(\tau,X_{\tau})\cdot \nabla F(\tau,X_{\tau}) 
|\nabla F \sigma(\tau,X_\tau)|
)
+o(\sqrt\Delta),
  \end{align*}
where $c_0$ is defined in \eqref{eq:c0}.
\end{thm}
Define now a smaller domain $\D^\Delta\subset \D$, which time-section is
given by $D^\Delta_t=\{x\in D_t: {\bf d}(x,\deD_t)> c_0 \sqrt \Delta
| \nabla F \sigma (t,
x)
|\}$, see Figure \ref{fig1}.
Introduce the exit time of the Euler scheme from this smaller domain: $\hat\tau^{\Delta}=\inf\{t_i> 0:X^\Delta_{t_i}\not \in D_{t_i}^\Delta\} \le \tau^\Delta$. The boundary correction procedure consists in simulating
\begin{align}
\label{ref_correc}
&  g(\hat\tau^{\Delta}\land T,
X^{\Delta}_{\hat\tau^{\Delta}\land T}
)
Z^{\Delta}_{\hat\tau^{\Delta}\land T}+\int_{0}^{\hat\tau^{\Delta}\land T} Z^{\Delta}_{\phi(s)} f(\phi(s),X^{\Delta}_{\phi(s)}) ds.
\end{align}
As above, we do not compute any projection on the boundary.
We denote the expectation of \eqref{ref_correc} by $u^\Delta(\D^\Delta)$.
One has:
\begin{thm}[Boundary correction]
\label{corr}
Under the assumptions of Theorem \ref{dev}, if we additionally suppose $\nabla F(.,.) | \nabla F \sigma(.,.)|$ is in ${\cal C}^{1,2}$, then one has:
$$u^\Delta(\D^\Delta)-u(\D)=o(\sqrt \Delta).$$
\end{thm}
The additional assumption is due to technical considerations to ensure that the modified domain $\D^\Delta$ is also of class $\H_2$. It is automatically fulfilled for domains of class ${\cal C}^3$ and $\sigma$ in ${\cal C}^{1,2}$.

\subsection{Proof of Theorems \ref{dev} and \ref{corr}}

\subsubsection{Error expansion}
By usual weak convergence arguments, Theorem \ref{dev} is a direct consequence of
Proposition \ref{tension} (tightness), Theorem \ref{loi_limite} (joint limit laws
associated to the overshoot) and Theorem \ref{FIRST_EXP} below.

\begin{thm}[First order approximation] Under the assumptions of
  Theorem \ref{dev}, one has
\label{FIRST_EXP}
  \begin{align*}
&u^\Delta(\D)-u(\D)=o(\sqrt {\Delta})+\\
&\E_x(\I_{\tau^{\Delta}\le T} Z_{\tau^\Delta}^\Delta (\nabla u-\nabla g)(\tau^\Delta,\pi_{\partial D_{\tau^\Delta}}(X^{\Delta}_{\tau^\Delta}))\cdot \nabla F(\tau^\Delta,X^{\Delta}_{\tau ^{\Delta}})
F^-(\tau^\Delta,X^{\Delta}_{\tau^\Delta})).
  \end{align*}
\end{thm}

\begin{rem}
In the above statement, we use projections on a non convex set, which needs a clarification. With the notation of Section \ref{time domains}, introduce $\tro:=\inf\{s> 0:X_s^\Delta\notin V_{D_s}(r_0)\}$. For $s\in[0,\tro]$ the projection
$\pi_{\bar {D}_s}(X_s^\Delta)$ is uniquely defined by
\begin{align}
\label{exp_proj}
\pi_{\bar {D}_s}(X^\Delta_{s})=X^\Delta_{s}+(\nabla F)^*(s,X^\Delta_{s})F^-(s,X^\Delta_{s}),
\end{align}
see Figure \ref{fig2}. Large deviation arguments (see Lemma \ref{contrprob} below)
also give $\P_x[\tro\le\tau^\Delta\le T]=O_{pol}(\Delta)$. Thus, in the following, for $s\ge \tau^{r_0} $, $\pi_{\bar D_{s}}(X_s^\Delta)$ and $\pi_{{ \partial D}_{s}}(X_s^\Delta) $ denote an arbitrary point on $\partial D_{s} $. This choice yields an exponentially small contribution in our estimates.
\psfrag{2PidD_t}{$\pi_{\partial D_t}(x)=y$}  
\psfrag{2x}{$x$}
\psfrag{2F(t,x)}{\hspace{-7mm}$F^-(t,x)$}
\psfrag{2nt}{$n(t,y)=[\nabla F]^*(t,x)$} 
\begin{figure}[h]  
\centerline{\epsfx{4cm}{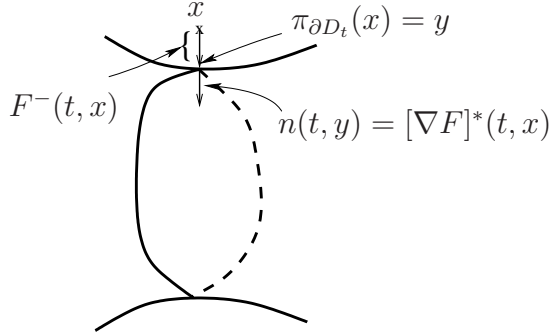}}  
\caption{Orthogonal projection $\pi_{\partial D_t}(x)$ of $x\notin D_t$ onto the boundary $\partial D_t$ and the related signed distance $F(t,x)$. Here $F(t,x)<0$ and ${\bf d}(x,\partial D_t)=|F(t,x)|=F^-(t,x)$.\label{fig2}} 
\end{figure} \end{rem}
{\sc Proof.}
Denote $e^\Delta:=u^\Delta(\D)-u(\D)$ the above error. Write now
\begin{align*}
e^\Delta=&\E_x[g(\tau^\Delta\land T,X_{\tau^\Delta\land T}^\Delta)Z_{\tau^ \Delta\land T}^\Delta-g(\tau^\Delta\land T,\pi_{\bar D_{\tau^\Delta\land T}}(X_{\tau^\Delta\land T}^\Delta))Z_{\tau^ \Delta\land T}^\Delta]\\
&+\big\{\E_x[g(\tau^\Delta\land T,\pi_{\bar D_{\tau^\Delta\land T}}(X_{\tau^\Delta\land T}^\Delta))Z_{\tau^ \Delta\land T}^\Delta+\int_{0}^{\tau^{\Delta}\land T} Z^{\Delta}_{\phi(s)} f(\phi(s),X^{\Delta}_{\phi(s)}) ds]\\
& \qquad -u(0,X_0^\Delta)\big\}\\
:=&e_1^\Delta+e_2^\Delta.
\end{align*}
We introduce here the projection for the error analysis. From \eqref{exp_proj} and Proposition \ref{tension}, a Taylor expansion yields
\begin{align}
e_1^\Delta=&-\E_x[\I_{\tau^\Delta\le T}Z_{\tau^\Delta}^\Delta\nabla g(\tau^\Delta, \pi_{\partial D_{\tau^\Delta}}(X_{\tau^\Delta}^\Delta))\cdot \nabla F(\tau^\Delta,X_{\tau^\Delta}^\Delta)F^-(\tau^\Delta,X_{\tau^\Delta}^\Delta) ]\nonumber\\
& +O(\Delta^{(1+\theta)/2}).\label{CTR_e1}
\end{align}

In the following, we write $U\=V$
(resp $U\<V$) when the equality between $U$ and $V$ holds in mean up
to a $O_{pol}(\Delta)$ (resp. $\E_x(U)\le\E_x(V)+O_{pol}(\Delta)$). We also use the notation $U=O(V)$ between two random variables $U$ and $V$ if for a constant $C$, one has $|U|\le C|V|$. Because $g(\tau^\Delta\land T,\pi_{\bar D_{\tau^\Delta\land T}}(X_{\tau^\Delta\wedge T}^\Delta))=u(\tau^\Delta\land T,\pi_{\bar D_{\tau^\Delta\land T}}(X_{\tau^\Delta\wedge T}^\Delta))$, we can write a telescopic summation:
\begin{align*}
  e_2^{\Delta}\=&  \bigl(\sum_{0\le t_i<\tau^\Delta\land T}
 u(t_{i+1},\pdtip(X^\Delta_{t_{i+1}}))Z^\Delta_{t_{i+1}}\\
&\qquad-u(t_{i},\pdti(X^\Delta_{t_{i}}))Z^\Delta_{t_i}
 + Z^\Delta_{t_i} f(t_i,X^\Delta_{t_{i}})\Delta\bigr)\I_{\tro>\tau^\Delta\land T}\\
 \=& \bigl(\sum_{0\le t_i< T}\1_{t_i<\tau^\Delta}\left[
 u(t_{i+1},\pdtip(X^\Delta_{t_{i+1}}))Z^\Delta_{t_{i+1}}\right.\\
&\qquad\left.-u(t_{i},X^\Delta_{t_{i}})Z^\Delta_{t_i} + Z^\Delta_{t_i} f(t_i,X^\Delta_{t_{i}})\Delta\right]\bigr)\I_{\tro>\tau^\Delta\land T}
\end{align*}
since for $t_i<\tau^\Delta$, $X^\Delta_{t_{i}}\in D_{t_i}$ and thus $\pdti(X^\Delta_{t_{i}})=X^\Delta_{t_{i}}$.  To proceed, the key idea is to introduce on the event 
$\{t_i<\tau^\Delta \}$, the partition $\{F(t_i,X_{t_i}^\Delta) \in (0,2\Delta ^{\frac12(1-\varepsilon)}]\}\cup\{F(t_i,X_{t_i}^\Delta)  > 2\Delta ^{\frac12(1-\varepsilon)} \} :=A_{t_i}^{\varepsilon}\cup (A_{t_i}^{\varepsilon})^C,\ \varepsilon>0$. This allows to split the cases
for which  $X_{t_i}^\Delta $  is close or not to the boundary $\partial D_{t_i} $. Lemma \ref{contrprob} ensures that $(X_{s}^\Delta)_{s\in[t_i,t_{i+1}]}$ stayed in $B(X_{t_i}^\Delta, \Delta^{\frac12(1-\varepsilon)})$ with a probability exponentially close to one. Then, on $(A_{t_i}^\varepsilon)^C $, the smoothness of the domain yields $\I_{(A_{t_i}^\varepsilon)^C} \P[X_{t_{i+1}}^\Delta \in D_{t_{i+1}}|\F_{t_i}]=1-O(\exp(-c\Delta^{-\varepsilon})) $, see Proposition \ref{int_shift} for a proof of this claim. On the other hand, on  $A_{t_i}^\varepsilon $, $X_{t_i}^\Delta $ is sufficiently close to the boundary to make the contribution of the overshoot at time $t_{i+1}$ significant for the error analysis. Write:
\begin{align}
  e_2^{\Delta}
 \=& \bigl(\sum_{0\le t_i< T}\1_{t_i<\tau^\Delta}\left\{\I_{A_{t_i}^\varepsilon}\left[
 u(t_{i+1},\pdtip(X^\Delta_{t_{i+1}}))Z^\Delta_{t_{i+1}}\right.\right. \nonumber \\
&\qquad\left.-u(t_{i},X^\Delta_{t_{i}})Z^\Delta_{t_i} + Z^\Delta_{t_i} f(t_i,X^\Delta_{t_{i}})\Delta\right]\nonumber \\
&+\I_{(A_{t_i}^\varepsilon)^C}\I_{\forall s\in[t_i,t_{i+1}],\ X_s^\Delta \in B(X_{t_i}^\Delta, \Delta^{\frac12(1-\varepsilon)})}\left[
 u(t_{i+1},X^\Delta_{t_{i+1}})Z^\Delta_{t_{i+1}}\right.\nonumber \\
&\qquad  \left.\left.-u(t_{i},X^\Delta_{t_{i}})Z^\Delta_{t_i} + Z^\Delta_{t_i} f(t_i,X^\Delta_{t_{i}})\Delta\right]\right\}\bigr)\I_{\tro>\tau^\Delta\land T}:=e_{21}^\Delta+e_{22}^\Delta. \label{DCOUP_EN_2}
\end{align}
Let us first deal with $e_{21}^\Delta$. In our framework, $u$ is $(1+\theta)/2 $-H\"older continuous in time and $\nabla u $ is $\theta $-H\"older continuous in space on a neighborhood of $\cD $. A Taylor expansion at order one and the equality \eqref{exp_proj} give
\begin{align*}
e_{21}^\Delta\=&\bigl(\sum_{0\le t_i< T}\1_{t_i<\tau^\Delta}\I_{A_{t_i}^\varepsilon}\left[Z_{t_i}^\Delta \nabla u(t_i,X_{t_i}^\Delta)\cdot\nabla F(t_{i+1},X_{t_{i+1}}^\Delta )F^-(t_{i+1},X_{t_{i+1}}^\Delta )\right. \\
&\left.+O(|F^-(t_{i+1},X_{t_{i+1}}^\Delta)|^{1+\theta})+O(|X^\Delta_{t_{i+1}}-X_{t_{i}}^\Delta|^{1+\theta})+O(\Delta^{\frac{1+\theta}{2}})\right] \bigr)\I_{\tro>\tau^\Delta\land T} \\
\=&\bigl(\I_{\tau^\Delta \le T} Z_{\tau^\Delta}^\Delta \nabla u(\tau^\Delta,X_{\tau^\Delta}^\Delta)\cdot \nabla F(\tau^\Delta,X_{\tau^\Delta}^\Delta)F^-(\tau^\Delta,X_{\tau^\Delta}^\Delta)\\ 
&+\sum_{0\le t_i< T}\1_{t_i<\tau^\Delta}\I_{A_{t_i}^\varepsilon}\bigl[ O(|F^-(t_{i+1},X_{t_{i+1}}^\Delta)|^{1+\theta})+O(|X^\Delta_{t_{i+1}}-X_{t_{i}}^\Delta|^{1+\theta})\\
&\qquad+O(|X_{t_{i+1}}^\Delta-X_{t_i}^\Delta|^\theta F^-(t_{i+1},X_{t_{i+1}}^\Delta))
+O(\Delta^{\frac{1+\theta}{2}})\bigr]\bigr) \I_{\tro>\tau^\Delta\land T}
\end{align*}
where we used once again Lemma \ref{contrprob} for the last equality. Standard arguments yield $\E[|X_{t_{i+1}}^\Delta-X_{t_i}^\Delta|^p|\F_{t_i}]=O(\Delta^{\frac p2}) $ for any $p>0$ and $ \E[|F^-(t_{i+1},X_{t_{i+1}}^\Delta)|^p|\F_{t_i}]=\E[|F^-(t_{i+1},X_{t_{i+1}}^\Delta)-F^-(t_i,X_{t_i}^\Delta)|^p|\F_{t_i}]=O(\Delta^{\frac p2})$ on $\{t_i<\tau^\Delta \}$. Thus, we can now rewrite
\begin{align*}
e_{21}^\Delta
\=&\bigl(\I_{\tau^\Delta \le T} Z_{\tau^\Delta}^\Delta \nabla u(\tau^\Delta,\pi_{\partial D_{\tau^\Delta}}(X_{\tau^\Delta}^\Delta))\cdot \nabla F(\tau^\Delta,X_{\tau^\Delta}^\Delta)F^-(\tau^\Delta,X_{\tau^\Delta}^\Delta)\bigr)\I_{\tro>\tau^\Delta\land T}\\
&+e_{211}^\Delta ,\\
e_{211}^\Delta\=&\bigl(\sum_{0\le t_i< T}\1_{t_i<\tau^\Delta}\I_{A_{t_i}^\varepsilon} O(\Delta^{\frac{1+\theta}{2}})
\bigr)\I_{\tro>\tau^\Delta\land T}.
\end{align*}
To handle $e_{211}^\Delta$ the idea is to use the occupation time formula and some sharp estimates concerning the local time of $(F(s,X_s^\Delta))_{s\le T\wedge  \tau^\Delta} $
 in a neighborhood of the boundary. We have
\begin{align*}
|e_{211}^\Delta|&\<C\Delta^{\frac{1+\theta }{2}}\biggl(\Delta^{-1}\bint{0}^{T\wedge
\tau^\Delta}\I_{F(\phi(t),X_{\phi(t)}^\Delta)\in
[0,2\Delta^{1/2(1-\varepsilon)}]
}dt\biggr)\I_{\tro>\tau^\Delta\land T} \\
&\<C\Delta^{\frac{1+\theta}{ 2}}\biggl(\Delta^{-1}\bint{0}^{T \wedge \tau^\Delta
}\I_{F(t,X_t^\Delta)\in [-\Delta^{1/2(1-\varepsilon)},3 \Delta^{1/2(1-\varepsilon)} ]
} dt \biggr)\I_{\tro>\tau^\Delta\land T}\\
&\<C\Delta^{\frac{1+\theta}{
2}}\biggl(\Delta^{-1}\bint{-\Delta^{1/2(1-\varepsilon)}}^{3\Delta^{1/2(1-\varepsilon)}}L_{T\wedge
\tau^\Delta}^y(F(.,X_{.}^\Delta)) dy\biggr)\I_{\tro>\tau^\Delta\land T},
\end{align*}
where we have used Lemma \ref{contrprob} at the second equality and the uniform ellipticity assumption for the last one.
Now an easy adaptation of the proof of Lemma 17 \cite{gobe:meno:spa:04} to our time-dependent domain framework gives
\begin{equation}
  \label{eq:tps:local}
  \E[L_{T\wedge \tau^\Delta}^y(F(.,X_{.}^\Delta))]\le C(|y|+\Delta^{\frac12}).
\end{equation}
Thus, one has $|e_{211}^\Delta| \< C\Delta^{\frac{1+\theta}{2}-\frac{\varepsilon}2}=o(\Delta^{\frac12})$ for $\varepsilon $ small enough. Hence, the above estimates and Lemma \ref{contrprob} give
\begin{align}
\label{CTR_E21}
e_{21}^\Delta&
\=\bigl(\I_{\tau^\Delta \le T} Z_{\tau^\Delta}^\Delta \nabla u(\tau^\Delta,\pi_{\partial D_{\tau^\Delta}}(X_{\tau^\Delta}^\Delta))\cdot \nabla F(\tau^\Delta,X_{\tau^\Delta}^\Delta)F^-(\tau^\Delta,X_{\tau^\Delta}^\Delta)\bigr)
+o(\Delta^{\frac12}).
\end{align}
Let us now turn to $e_{22}^\Delta $. 
If $g\in \H_{3+\theta}$ (which implies $u\in \H_{3+\theta}$ in view of Proposition \ref{prop:feynman:kac:stopped}), the term $e_{22}^\Delta$ can be handled with somehow standard techniques. Namely Taylor like expansions in the spirit of Talay and Tubaro
\cite{tala:tuba:90}. For simplicity we handle $e_{22}^\Delta $ under the previous smoothness assumption on $g$ and $u$. The proof under weaker assumptions ($g\in \H_{1+\theta}$), that involves sharp estimates on possibly exploding derivatives of $u$ near the boundary, is postponed to the Appendix. We recall that 
\begin{align*}
e_{22}^\Delta&
\=\bigl(\sum_{0\le t_i< T}\1_{t_i<\tau^\Delta}\I_{(A_{t_i}^\varepsilon)^C}\I_{\forall s\in[t_i,t_{i+1}],\ X_s^\Delta \in B(X_{t_i}^\Delta, \Delta^{\frac12(1-\varepsilon)})}\left[
 u(t_{i+1},X^\Delta_{t_{i+1}})Z^\Delta_{t_{i+1}}\right.\nonumber \\
&\qquad \qquad\left.-u(t_{i},X^\Delta_{t_{i}})Z^\Delta_{t_i} + Z^\Delta_{t_i} f(t_i,X^\Delta_{t_{i}})\Delta\right]\bigr)\I_{\tro>\tau^\Delta\land T} 
\end{align*}
For all $(s,y)\in \cD $ introduce the operators $L_{s,y}:{\cal C}^{1,2}(\cD)\rightarrow {\cal C}(\cD), \ \varphi\mapsto ((t,x)\mapsto L_{s,y}\varphi(t,x)=\nabla \varphi(t,x) b(s,y)+\frac 12 \Tr[H\varphi(t,x) [\sigma\sigma^*](s,y)]  ) $.
Recalling that $\partial_t u(t_i,X_{t_i}^\Delta)+L_{t_i,X_{t_i}^\Delta} u(t_i,X_{t_i}^\Delta)-ku(t_i,X_{t_i}^\Delta)+f(t_i,X_{t_i}^\Delta)=0 $, It\^o's formula gives
\begin{align}
e_{22}^\Delta
\=&\bigl(\sum_{0\le t_i< T}\1_{t_i<\tau^\Delta}\I_{(A_{t_i}^\varepsilon)^C}\I_{\forall s\in[t_i,t_{i+1}],\ X_s^\Delta \in B(X_{t_i}^\Delta, \Delta^{\frac12(1-\varepsilon)})}\bigl[\nonumber \\
&\bint{t_i}^{t_{i+1}} (Z_s^\Delta-Z_{t_i}^\Delta)(\partial_s+L_{t_i,X_{t_i}^\Delta}-
 k(t_i,X_{t_i}^\Delta))u(s,X_s^\Delta)ds\nonumber\\
 &+Z_{t_i}^\Delta\bint{t_i}^{t_{i+1}} \bigl[\bigl(\partial_s+L_{t_i,X_{t_i}^\Delta}-k(t_i,X_{t_i}^\Delta)\bigr)u(s,X_s^\Delta)\nonumber\\
&\hspace{2cm} -\bigl(\partial_s+L_{t_i,X_{t_i}^\Delta}-k(t_i,X_{t_i}^\Delta)\bigr)u(t_i,X_{t_i}^\Delta)) \bigl]ds\nonumber\\
 &+M_{t_i,t_{i+1}}\bigr]\bigr)\I_{\tro>\tau^\Delta\land T},\label{PREAL_DEP_DER_VARIEES}
 \end{align}
where for all $v\in[t_i,t_{i+1}] $, $M_{t_i,v}:= \bint{t_i}^{v} Z_s^\Delta \nabla u(s,X_s^\Delta)\sigma(t_i,X_{t_i}^\Delta) dW_s$ is a square-integrable martingale term. Note that in this definition, in whole generality, $M_{t_i,v} $ is not stopped at the exit time $\tau_{t_i}:=\inf\{s\ge t_i: X_s^\Delta \not\in D_s \} $.  If $\tau_{t_i}\le t_{i+1} $ (which happens with exponentially small probability on $(A_{t_i}^\varepsilon)^C $), the term $\nabla u(s,X_s^\Delta),\ s\in[\tau_{t_i},t_{i+1}] $ in $M_{t_i,t_{i+1}} $ has to be understood as the smooth extension of $\nabla u $ to the whole space. In particular this extension remains bounded. 
Now, we derive from Lemma \ref{contrprob}
\begin{align*}
\E_x[\bigl(\sum_{0\le t_i< T}\1_{t_i<\tau^\Delta}\I_{(A_{t_i}^\varepsilon)^C}\I_{\forall s\in[t_i,t_{i+1}],\ X_s^\Delta \in B(X_{t_i}^\Delta, \Delta^{\frac12(1-\varepsilon)})} M_{t_i,t_{i+1}}\bigr]\bigr)\I_{\tro>\tau^\Delta\land T}]\\
=\E_x[\sum_{0\le t_i< T}\1_{t_i<\tau^\Delta}\I_{(A_{t_i}^\varepsilon)^C}M_{t_i,t_{i+1}}] +O_{pol}(\Delta)=O_{pol}(\Delta).
\end{align*}

We can thus neglect the contribution of the martingale terms in \eqref{PREAL_DEP_DER_VARIEES}. We now develop the other quantities in \eqref{PREAL_DEP_DER_VARIEES} with Taylor integral formulas to derive
\begin{align}
&\bint{t_i}^{t_{i+1}} (Z_s^\Delta-Z_{t_i}^\Delta)(\partial_s+L_{t_i,X_{t_i}^\Delta}- k(t_i,X_{t_i}^\Delta))u(s,X_s^\Delta)ds\nonumber\\
 =&O\bigl(\Delta^2(|u|_\infty+|\nabla u|_\infty+|\partial_t u|_\infty+|D^2 u|_\infty )\bigr),\nonumber\\
 &\bint{t_i}^{t_{i+1}} (\partial_tu(s,X_s^\Delta)-\partial_tu(t_i,X_{t_i}^\Delta)) ds\nonumber\\
=&\bint{t_i}^{t_{i+1} } \nabla \partial_t u(t_i,X_{t_i}^\Delta) \sigma(t_i,X_{t_i}^\Delta) (W_s-W_{t_i})ds \nonumber\\
 &+O\bigl(\Delta^{1+\frac{1+\theta}{2}}[\partial_t u]_{t,\frac{1+\theta}{2}}+\Delta^2|\nabla\partial_t  u|_\infty+\Delta \sup_{s\in[t_i,t_{i+1}]}|X_s^\Delta-X_{t_i}^\Delta|^{1+\theta}[\nabla\partial_t   u]_{x,\theta} \bigr),\nonumber
\end{align}
\begin{align}
& \bint{t_i}^{t_{i+1}} (L_{t_i,X_{t_i}^\Delta}u(s,X_s^\Delta)-L_{t_i,X_{t_i}^\Delta}u(t_i,X_{t_i}^\Delta)) ds\nonumber\\
=&\bint{t_i}^{t_{i+1}}\langle H_u(t_i,X_{t_i}^\Delta) \sigma(t_i,X_{t_i}^\Delta) (W_s-W_{t_i}),b(t_i,X_{t_i}^\Delta)\rangle  ds \nonumber\\
 &+\frac 12 \bint{t_i}^{t_{i+1}}  \Tr\bigl( (D^3 u(t_i,X_{t_i}^\Delta) \sigma(t_i,X_{t_i}^\Delta)(W_s-W_{t_i})) \cdot a(t_i,X_{t_i}^\Delta)\bigr) ds\nonumber\\
 &+O\bigl(\Delta^2 \{ |D^2u|_\infty+|D^3u|_\infty+|\partial_t \nabla u|_\infty\}+\Delta^{1+\frac{1+\theta}2}[D^2 u]_{t,\frac{1+\theta}2}\nonumber\\
 &+\Delta |D^3 u|_\infty\sup_{s\in[t_i,t_{i+1}]} |X_s^\Delta-X_{t_i}^\Delta|^2+\Delta \sup_{s\in[t_i,t_{i+1}]}|X_{s}^\Delta-X_{t_i}^\Delta|^{1+\theta}[D^3u]_{x,\theta} \bigr),\nonumber\\
&k(t_i,X_{t_i}^\Delta)\bint{t_i}^{t_{i+1}}(u(s,X_s^\Delta)-u(t_i,X_{t_i}^\Delta)) ds\nonumber\\
=&k(t_i,X_{t_i}^\Delta)\bint{t_i}^{t_{i+1} } \nabla   u(t_i,X_{t_i}^\Delta) \sigma(t_i,X_{t_i}^\Delta) (W_s-W_{t_i})ds\nonumber \\
 &+O\bigl(\Delta^{2}(|\partial_t  u|_\infty+|\nabla  u|_\infty)
 +\Delta |D^2 u|_\infty\sup_{s\in[t_i,t_{i+1}]}|X_s^\Delta-X_{t_i}^\Delta|^{2}
 \bigr), \label{THE_DETAILS}
\end{align}
where $[\cdot]_{t,\alpha}, [\cdot]_{x,\alpha},\alpha \in(0,1]$ denote respectively the H\"older norms of order $\alpha $ in time and space (see Chapter IV Section 1 p. 46 in \cite{lieb:96} for a precise definition).

Hence, bringing together our estimates and exploiting the relations between the spatial and time derivatives for $u$ (through the PDE), from (\ref{PREAL_DEP_DER_VARIEES}) and (\ref{THE_DETAILS}) we derive
 \begin{align}
e_{22}^\Delta
\= & \bigl(\sum_{0\le t_i< T}\1_{t_i<\tau^\Delta}\I_{(A_{t_i}^\varepsilon)^C}\I_{\forall s\in[t_i,t_{i+1}],\ X_s^\Delta \in B(X_{t_i}^\Delta, \Delta^{\frac12(1-\varepsilon)})}\bigl[\nonumber\\
& \quad O\bigl (\Delta^2\{1+|u|_\infty+|\nabla u|_\infty+|D^2 u|_\infty+|D^3 u|_\infty\}\bigr)\nonumber\\
& +O\bigl (\Delta^{1+\frac{1+\theta}{2}}\{1+|u|_\infty+|\nabla u|_\infty+|D^2 u|_\infty+|D^3 u|_\infty+[D^2 u]_{t,\frac{1+\theta}{2}}\}\bigr)\nonumber\\
&+O\bigl (\Delta \sup_{s\in[t_i,t_{i+1}]}|X_s^\Delta-X_{t_i}^\Delta|^{1+\theta}\{1+|u|_\infty+|\nabla u|_\infty+|D^2 u|_\infty+|D^3 u|_\infty+[D^3 u]_{x,\theta}\}\bigr)\nonumber\\
&+O\bigl (\Delta \sup_{s\in[t_i,t_{i+1}]}|X_s^\Delta-X_{t_i}^\Delta|^{2}\{|D^2 u|_\infty+|D^3 u|_\infty\}\bigr)\nonumber\\
&+\bar M_{t_i,t_{i+1}}\bigr]\bigr)\I_{\tro>\tau^\Delta\land T},\label{DEP_DER_VARIEES}
\end{align}
 where $\bar M_{t_i,t_{i+1}} $ denotes the sum of the terms involving the Brownian increment $(W_s-W_{t_i})_{s\in[t_i,t_{i+1}]} $ in the above equations \eqref{THE_DETAILS}.
 Under our current assumption, i.e. $u\in \H_{3+\theta} $, all the norms appearing in \eqref{DEP_DER_VARIEES} and all the derivatives appearing in the $(\bar M_{t_i,t_ {i+1} })_{0\le t_i<T} $ are bounded. Hence,
 \begin{align}
\1_{t_i<\tau^\Delta}&\I_{(A_{t_i}^\varepsilon)^C}\E[\I_{\forall s\in[t_i,t_{i+1}],\ X_s^\Delta \in B(X_{t_i}^\Delta, \Delta^{\frac12(1-\varepsilon)})} \bar M_{t_i,t_{i+1}}\I_{\tro>\tau^\Delta\land T}|\F_{t_i}]\nonumber\\
\qquad =&\1_{t_i<\tau^\Delta}\I_{(A_{t_i}^\varepsilon)^C}\E[\bar M_{t_i,t_{i+1}}|\F_{t_i}] +O_{pol}(\Delta)=O_{pol}(\Delta),\label{CTR_E22_0}\\
|e_{22}^\Delta|\<&  C\Delta^{\frac{1+\theta}2}.\label{CTR_E22}
\end{align}
Plug \eqref{CTR_E21} and \eqref{CTR_E22} into \eqref{DCOUP_EN_2}. The statement is derived from \eqref{CTR_e1} and \eqref{DCOUP_EN_2}. 
We specify in the Appendix how to complete the proof from a sharper version of \eqref{DEP_DER_VARIEES} deriving from \eqref{PREAL_DEP_DER_VARIEES}, when $g\in \H_{1+\theta}$. \qed

\subsubsection{Boundary Correction}
\label{par:boundary:correction}
One has
\begin{align}
  u^\Delta(\D^\Delta)-u(\D)=[u^\Delta(\D^\Delta)-u(\D^\Delta)]+[u(\D^\Delta)-u(\D)].
\label{cont:1}
\end{align}
\begin{enumerate}
\item The first contribution in \eqref{cont:1} has been previously analysed in Theorem \ref{dev}, except that the domain $\D^\Delta$ depends on $\Delta$. We can show that it is equal to $c_0 \sqrt \Delta \E(\1_{\tau\le T} Z_\tau  (\nabla
u-\nabla g)(\tau,X_{\tau})\cdot \nabla F(\tau,X_{\tau}) | \nabla
F \sigma(\tau,X_\tau)|)+o(\sqrt\Delta)$.\\ We briefly sketch the proof of this assertion, which is done in two steps. For this, set $\hat u^\Delta=u(\D^\Delta)$ for the solution of the PDE in the domain $\D^\Delta$.
\begin{itemize}
\item {\bf Step 1.} It is well known that all PDE estimates depend only on bounds on the derivatives of the level set functions $(\varphi_0)$ arising in the definition of the time-dependent domains (see section \ref{time domains}), and on the bounds on the derivatives of data $g$, $f$ and $k$. Hence, since $\D^\Delta$ is a small perturbation of class $\H_2$ (because $\nabla F|\nabla F\sigma|$ has this regularity) of the domain $\D$ of class $\H_2$, all PDE estimates on $\hat u^\Delta$ remain locally uniform w.r.t. $\Delta$. In addition, $\hat u^\Delta$ and its gradient converge uniformly to $u$ and $\nabla u$. This argumentation allows us to state that the first order approximation theorem holds:
  \begin{align*}
&u^\Delta(\D^\Delta)-u(\D^\Delta)=o(\sqrt {\Delta})+\\
&\E_x(\I_{\hat\tau^{\Delta}\le T} Z_{\hat \tau^\Delta}^\Delta (\nabla u-\nabla g)( \hat \tau^\Delta,\pi_{\partial D^\Delta_{\hat \tau^\Delta}}(X^{\Delta}_{\hat \tau^\Delta}))\cdot \nabla \hat F^\Delta(\hat \tau^\Delta,X^{\Delta}_{\hat \tau ^{\Delta}})
[\hat F^\Delta]^-(\hat \tau^\Delta,X^{\Delta}_{\hat \tau^\Delta})),
  \end{align*}
where $\hat F^\Delta$ and $\hat \tau^\Delta$ are respectively the signed distance to the side of $\D^\Delta$ and the related discrete exit time.
\item {\bf Step 2.} The second step is to prove that the analogous version of Theorem \ref{loi_limite} holds, with $\hat \tau^\Delta$ instead of $\tau^\Delta$. Actually, a careful reading of its proof shows that it is indeed the case, without modification.
\end{itemize}
\item Finally, the last term in \eqref{cont:1} is related to the sensitivity of a Dirichlet problem with respect to the domain. By an application of Theorem 2.2 in \cite{cost:elka:gobe:06} with $\Theta(t,x)=-c_0 \nabla F(t,x) | \nabla F \sigma(t,x)|$ (in ${\cal C}^{1,2}$), one gets that this contribution equals
$$-c_0 \sqrt \Delta \E(\1_{\tau\le T}Z_\tau (\nabla u-\nabla
g)(\tau,X_{\tau})\cdot \nabla F(\tau,X_{\tau}) | \nabla
F\sigma(\tau,X_\tau)|)+o(\sqrt\Delta).$$
\end{enumerate}
This proves that the new procedure has an error $o(\sqrt \Delta)$.\qed

\mysection{Technical results concerning the overshoot}

\label{proof_res_os}
This section is devoted to the proof of Theorem
\ref{loi_limite}.
We first state some useful auxiliary results.

\begin{lem} [Bernstein's inequality]
\label{contrprob} Assume \A{{A$_{\theta}$-1}}. Consider two stopping
times $S,S'$ upper bounded by $T$ with $0\le S'-S\leq \Theta \leq T
$. Then for any $p\geq 1$, there are some constants $c>0$ and $C:=C($\A{{A$_{\theta}$-1}} , $T)$,
such that for any $\eta\geq 0$, one has $a.s$:
\begin{align*}
  \P[\sup_{t\in [S,S']}|X_t^\Delta-X_S^\Delta|\geq \eta\;\big|\;\F_S] \leq & C\exp\left( -c\bfrac{\eta^2}{\Theta}\right),\\
\E[\sup_{t\in[S,S']}| X_t^\Delta- X_S^\Delta|^p\;\big|\;\F_S]\leq & C
\Theta^{p/2}.
\end{align*}
\end{lem}
For a proof of the first inequality we refer to Chapter 3, \S 3 in
\cite{revu:yor:99}. The last inequality easily follows from the
first one or from the BDG inequalities.

\begin{lem} [Convergence of exit time]
\label{conv:proba:tau:xtau} Assume \A{{A$_{\theta}^{'}$}} and that the domain is of class $\H_2$. The following convergences hold in probability:
\begin{enumerate}
\item $\lim_{\Delta\rightarrow 0}\tau^\Delta\land T=\tau\land T$;
\item $\lim_{\Delta\rightarrow 0}X^\Delta_{\tau^\Delta\land T}=X_{\tau\land T}$;
\item $\lim_{\Delta\rightarrow 0}\sup_{t\le T}|X^\Delta_{\phi(t)}-X_t|=0$.
\end{enumerate}
\end{lem} The proof of the first two assertions in the case of space-time domain is analogous to the case of cylindrical domain (see \cite{gobe:mair:05}) and thus left to the reader. The last convergence is standard.

The following results are key tools to prove Theorem 
\ref{loi_limite}. A similar version is proved in \cite{sieg:79}, but here, we additionally prove the uniform convergence.
\begin{lem}{\bf (Asymptotic independence of the overshoot and the discrete exit time).} \label{lemme_siegmund}
Let $W$ be a standard one dimensional BM. Put $x>0$ and
consider the domain $\cD:=]0,T[\times ]-\infty,x[$. With the notation of Section
\ref{main_results}, for any $\varepsilon>0$ we have
\begin{align}
\label{indas} &\lim_{\Delta\longrightarrow 0} \sup_{t\in[0,T],y\ge 0,x\ge \Delta^{1/2-\varepsilon}}\left|\P_0[ \tau^\Delta\le
t, (W_{\tau^\Delta}-x)\le y\sqrt \Delta]-\P_0[\tau\leq  t]H(y)\right|=0.
\end{align}
\end{lem}

If the Euler scheme starts close to the boundary at a small distance $d$, its discrete exit likely occurs after a time roughly equal to $d^2$. This feature is quantified in the above lemma.
\begin{lem}
\label{CTR_Hitting} Assume \A{$\mathbf{A}_\theta^{'} $}, and that
the domain is of class $\H_2$. Let $0<\beta<\alpha<1/2 $. For all
$\eta>0$, there exists $C:=C_\eta>0$ s.t. for $\Delta$ small enough,
$\forall s\in \Delta \mathbb N\cap[0,T]$ and $\forall x\in V_{\partial D_s}(\Delta^\alpha)\cap D_s$, one has
\begin{align*}
\P[\tau^\Delta\land T \ge \Delta^{2\beta }|X^\Delta_s=x] \le C(\Delta ^{\alpha-\beta-\eta}+\Delta^\beta),
\end{align*}
where $\tau^\Delta:=\inf\{t_i> s:X_{t_i}^\Delta\notin
D_{t_i}\}$.
\end{lem}

\begin{lem}
\label{CTR_RESTES} Assume \A{$\mathbf{A}_\theta^{'} $}, and that the
domain is of class $\H_2$. There exists
$ C>0$, such that $\forall s\in \Delta \mathbb N\cap[0,T]$, $\forall x\in D_s$, $\forall t\in[s,T]$
and $\forall b\ge a\ge 0$, one has
\begin{align*}
\P[\tau^\Delta\le
  t,\Delta^{-1/2}F^-(\tau^\Delta,X_{\tau^\Delta}^\Delta) \in [a,b]
|X^\Delta_s=x]\le& C\bigl((b-a)+\Delta^{1/4}\bigr)
\end{align*}
where $\tau^\Delta$ is shifted as in the previous lemma.
\end{lem}
The proof of these three lemmas is postponed to Section
\ref{le_inter}.

We mention that if $\sigma\sigma^* $ is uniformly elliptic, 
Lemma \ref{CTR_RESTES} is valid without the $\Delta^{1/4}$ (see the proof for details). In that case, it means that the law of the renormalized overshoot is absolutely continuous w.r.t. the Lebesgue measure on $\R^+$, with a bounded density. This is also true at the limit, in view of Theorem \ref{loi_limite}. 

\subsection{Proof of Theorem \ref{loi_limite}}
\label{ref_pre_ll} 
\begin{figure}[h]  
\centerline{\epsfx{6cm}{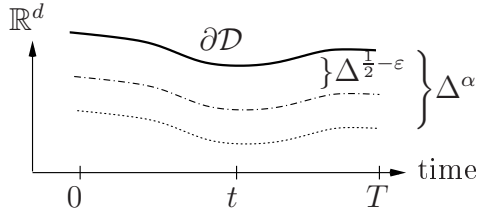}}  
\caption{The two localization neighborhoods with $\alpha<\frac 12-\varepsilon$.\label{fig5}} 
\end{figure}

Consider first the case $\cD=]0,T[\times D$ where $D$ is a half space. The theorem in the case of BM is then a direct consequence of Lemma \ref{lemme_siegmund}. Now to deal with the Euler scheme, we introduce
a first neighborhood whose distance to the boundary goes to 0 with
$\Delta$ at a speed lower than $\Delta^{1/2}$ (below, the speed is tuned by a parameter $\alpha$, see Figure \ref{fig5}). The characteristic exit time for a starting point in this neighborhood is short (Lemma \ref{CTR_Hitting}), thus the diffusion coefficients are somehow constant and we are almost in the BM framework. Also, a second localization
w.r.t. to the hitting time of this neighborhood guarantees that up
to a rescaling we are far enough from the boundary to apply the renewal
arguments needed for the asymptotic law of the overshoot (this is tuned by another parameter $\varepsilon$, see Figure \ref{fig5}).\\
For a more general time-space domain of class $\H_2$ two additional tools are used: a time-space change of chart and a local half space approximation of the domain by some tangent hyperplane.\\
For notational convenience, we  assume from now on that the time-section domains 
$(D_t)_{t\in[0,T]}$ are convex so that $\pi_{\partial D_t}$ is always
uniquely defined on $D_t^c $. To handle the case of general $\H_2$ domains,
an additional localization procedure similar to the one of Theorem \ref{FIRST_EXP} is needed. We leave it to the reader.\\
For the sake of clarity, we also assume $k\equiv 0$ ($Z\equiv 1$). This is an easy simplification since owing to Lemma \ref{conv:proba:tau:xtau}, $Z^\Delta_{\tau^\Delta\land T}$ converges to $Z_{\tau\land T}$ in $L_1$.

\textbf{Step 1: preliminary localization.} For $\alpha <1/2$ specified later on,
define $\tha:=\inf\{t_i> 0:F(t_i,X_{t_i}^{\Delta})\le
\Delta^\alpha\}\le \tau^\Delta$. We aim at studying the convergence of
$$\varPsi_\Delta(t,x,y):=\E_x[\I_{\tau^\Delta\le
t,F^-(\tau^\Delta,X_{\tau^\Delta}^{\Delta})\ge y\sqrt
\Delta}\varphi(X_{\tau^\Delta}^\Delta)]$$
and for this, we define for all $0\le s\le t 
<T$ ($s\in \Delta \mathbb N$), $(\tilde x,y)\in\R^d\times \R^+$
\begin{align*}
\Psi_\Delta(s,t,\tilde x,y):=&\P[{\tau^\Delta\le
t,F^-(\tau^\Delta,X_{\tau^\Delta}^{\Delta})\ge y\sqrt \Delta}
|X_s^\Delta=\tilde x],\\
A(t,\alpha,\varepsilon):=&\{\tha<\tau^\Delta,\tha< t,F(\tha,X_{\tha}^{\Delta})\ge  \Delta^{1/2-\varepsilon}\}.
\end{align*}
Here, $\varepsilon$ is a fixed parameter in $]0,1/2[$, such that $\alpha<1/2-\varepsilon$ (take $\varepsilon=(\alpha+1/2)/2$ for instance).\\
In the definition of $\Psi_\Delta$, $\tau^\Delta$ has to be understood as the shifted exit time $\inf\{t_i> s:X_{t_i}^\Delta\notin D_{t_i}\}$.
By Lemma \ref{contrprob}, $\P_x[\tau^\Delta=\tau_{\Delta^\alpha}\le
t]+\P_x[\tha< t,F(\tha,X_{\tha}^{\Delta})<
  \Delta^{1/2-\varepsilon}]=O_{pol}(\Delta)$ using $\alpha< 1/2-\varepsilon$. Hence,
\begin{align*}
\varPsi_\Delta(t,x,y)
=& \E_x[\I_{A(t,\alpha,\varepsilon),F^-(\tau^\Delta,X_{\tau^\Delta}^{\Delta})\ge y\sqrt \Delta}\varphi(X_{\tau^\Delta}^\Delta)\1_{\tau^\Delta\le t}]+O_{pol}(\Delta)\nonumber\\
=&\E_x[\I_{A(t,\alpha,\varepsilon),F^-(\tau^\Delta,X_{\tau^\Delta}^{\Delta})\ge y\sqrt \Delta}(\varphi(X_{\tau^\Delta}^\Delta)-\varphi(X_{\tau_{\Delta^\alpha}}^\Delta))\1_{\tau^\Delta\le t}]\nonumber\\
&+\E_x[\I_{A(t,\alpha,\varepsilon)} \varphi(X_{\tau_{\Delta^\alpha}}^\Delta)\Psi_\Delta(\tha,t,X_{\tha}^\Delta,y)]+O_{pol}(\Delta).
\end{align*}
The first term in the right hand side above converges to 0, using the convergence in probability of $|X^\Delta_{\tau^\Delta\land T}-X_{\tau_{\Delta^\alpha}\land T}^\Delta|$ to 0 (analogously to Lemma \ref{conv:proba:tau:xtau}). This gives
\begin{equation}
  \label{first_ctr}
  \varPsi_\Delta(t,x,y)
= \E_x[\I_{A(t,\alpha,\varepsilon)} \varphi(X_{\tau_{\Delta^\alpha}}^\Delta)\Psi_\Delta(\tha,t,X_{\tha}^\Delta,y)]+o(1).
\end{equation}

\psfrag{5t}{$t$}  
\psfrag{50}{$0$}  
\psfrag{5T}{$T$} 
\psfrag{5Rd}{$\R^d$} 
\psfrag{5time}{time} 
\psfrag{5partialD}{$\partial \D$}  
\psfrag{5eps}{$\}\Delta^{\frac 12-\varepsilon}$}
\psfrag{5alpha}{$\bigg\}\Delta^{\alpha}$}
Let us comment again these two localisations. That with $\Delta^\alpha$ enables us to freeze the coefficients of the Euler scheme, because the exit time is likely close to the initial time. That with $\Delta^{1/2-\varepsilon}$ ensures that it starts far enough from the boundary to induce the limiting behavior of the overshoot. This right balance regarding the distance of the initial point to the boundary is crucial. The final choice of $\alpha$ (and thus $\varepsilon$) depends on the regularity $\theta$ of the coefficients $b$ and $\sigma$.\\
Now, it remains to study the convergence of $\Psi_\Delta(.)$.

\textbf{Step 2: diffusion with frozen coefficients.} Denote
$\tha:=\tilde s,\ X_\tha^\Delta:=\tilde x$. Conditionally to
$\F_{\tilde s}$, introduce now the one dimensional process
$(Y_s)_{s\ge\tilde s}$, $Y_s=F(\tilde
s,\tilde x)
+(\nabla F
\sigma)(\tilde s,\tilde x)(W_s-W_{\tilde
s})$. Note that we do not take into account the drift part in the
frozen process. From the next localization procedure, it yields a
negligible term. Since $Y$ has constant coefficients, we apply below
Lemma \ref{lemme_siegmund} to handle the overshoot of $Y$ w.r.t.
$\R^{+*}$. Define $\tau^{\Delta,Y}:=\inf\{t_i>\tilde s: Y_{t_i}
\le 0\}  $ and rewrite
\begin{align}
\label{ctr_psi_c} \Psi_\Delta(\tilde s,t,\tilde
x,y)&:=\Psi_\Delta^C(\tilde s,t,\tilde x,y
  )+R_\Delta(\tilde s,
  t,\tilde x,y),\\
   \Psi_\Delta^C(\tilde s,t,\tilde x,y)&:=\P_{\tilde s,\tilde
  x}[\tau^{\Delta,Y}\le t,(Y_{\tau^{\Delta,Y}})^-\ge y\sqrt \Delta].\nonumber
\end{align}
From \A{$\mathbf{A}_\theta^{'} $-\textbf{2'}} that guarantees that
$Y$ has a non degenerate variance and Lemma \ref{lemme_siegmund},
one gets
\begin{align*}
\sup_{(\tilde s,\tilde x)\in \cal A^{\alpha,\varepsilon}
} |\Psi_\Delta^C(\tilde s,t,\tilde x, y)&-\P_{\tilde s,\tilde
x}[\tau^{\Delta,Y}\le
  t](1-H(y/
|(\nabla F \sigma )(\tilde s,\tilde x )  
|))|
\underset{\Delta \rightarrow 0}{\longrightarrow} 0,
\end{align*}
where ${\cal A}^{\alpha,\varepsilon}:=\{(t,x):0\le t \le T, x\in V_{\partial D_t(\Delta^\alpha)}\backslash V_{\partial D_t}(\Delta^{1/2-\varepsilon})\}.$
Plug now this identity in
\eqref{ctr_psi_c} to obtain with the same uniformity
\begin{align}
\label{CTR_PSI} \Psi_\Delta(\tilde s,t,\tilde x,y)=& \P_{\tilde
s,\tilde x}[\tau^{\Delta,Y}\le
  t](1-H(y/
|( \nabla F \sigma)(\tilde s,
\tilde x
) 
|))
+ R_\Delta(\tilde s,t,\tilde x,y) +o(1).
\end{align}
\textbf{Step 3: control of the rests.} We now show that
$R_\Delta(\tilde s,t,\tilde x,y)=o(1)$ where the rest is still
uniform for $(\tilde s,\tilde x)\in
{\cal A}^{\alpha,\varepsilon}
$. This part is long and technical. First, decomposing the space using the events $\{\tau^\Delta=\tau^{\Delta,Y}\}$, $\{F^-(\tau^\Delta,X_{\tau^\Delta}^{\Delta} )\ge y\sqrt \Delta\}$, $\{(Y_{\tau^{\Delta,Y}})^-\ge y\sqrt \Delta\}$ and their complementary events, write: 
\begin{align}
\label{reste_1} |R_\Delta &(\tilde s,t,\tilde x,y)| \le R_\Delta^1(\tilde s,
  t,\tilde x)\nonumber\\
&+\P_{\tilde s,\tilde
x}[\tau^\Delta\le
  t,F^-(\tau^\Delta,X_{\tau^\Delta}^{\Delta} )\ge y\sqrt
  \Delta,(Y_{\tau^{\Delta,Y}})^-< y\sqrt \Delta,\tau^\Delta=\tau^{\Delta,Y}]\nonumber \\
& +\P_{\tilde s,\tilde
x}[\tau^\Delta\le
  t,F^-(\tau^\Delta,X_{\tau^\Delta}^{\Delta} )< y\sqrt
  \Delta,(Y_{\tau^{\Delta,Y}})^-\ge y\sqrt \Delta,\tau^\Delta=\tau^{\Delta,Y}]
\end{align}
with $R_\Delta^1(\tilde s,t,\tilde x)\le \P_{\tilde s,\tilde
x}[\tau^\Delta\le t,
  \tau^\Delta\neq \tau^{\Delta,Y}]+\P_{\tilde s,\tilde x}[\tau^{\Delta,Y}\le t,
  \tau^\Delta\neq \tau^{\Delta,Y}]:=(R_\Delta^{11}+R_\Delta^{12})(\tilde s,t,\tilde x)$. Let $y_\Delta $ be a given positive
  function of the time-step s.t. $y_\Delta\underset{\Delta\rightarrow 0}{\rightarrow} 0 $ specified
  later on.\\
On the event $\{\tau^\Delta=\tau^{\Delta,Y},|Y_{\tau^{\Delta,Y}}
 -F(\tau^{\Delta,Y},X_{\tau^{\Delta,Y}}^{\Delta})|\le y_\Delta\sqrt \Delta \}$, the conditions $F^-(\tau^\Delta,X_{\tau^\Delta}^{\Delta})\ge y\sqrt
  \Delta $ and $(Y_{\tau^{\Delta,Y}})^-<y\sqrt {\Delta} $ imply $\Delta^{-1/2}(Y_{\tau^{\Delta,Y}})^-\in [y-y_\Delta, y)$. Similarly, $(Y_{\tau^{\Delta,Y}})^-\ge y\sqrt {\Delta} $ and $F^-(\tau^\Delta,X_{\tau^\Delta}^{\Delta})< y\sqrt
  \Delta $ imply $\Delta^{-1/2}(Y_{\tau^{\Delta,Y}})^-$ $\in [y,y+y_\Delta) $. Hence, by setting
\begin{align*}
R_\Delta^2(\tilde s,t,\tilde x)&:= 2 \P_{\tilde s,\tilde
  x}[\tau^{\Delta,Y}\le t,\tau^\Delta=\tau^{\Delta,Y},|Y_{\tau^{\Delta,Y}}
 -F(\tau^{\Delta,Y},X_{\tau^{\Delta,Y}}^{\Delta})|> y_\Delta \sqrt \Delta    ],\\
R_\Delta^3(\tilde s,t,\tilde x,y)&:=\P_{\tilde s,\tilde
x}[\tau^{\Delta,Y}\le t,
 \Delta^{-1/2}(Y_{\tau^{\Delta,Y}} )^-\in[y-y_\Delta,y+y_{\Delta}),\tau^\Delta=\tau^{\Delta,Y}],
\end{align*}
we obtain
$R_\Delta(\tilde s,t,\tilde x,y)|\le (R_\Delta^1+R_\Delta^2)(\tilde s,t,\tilde x) +R_\Delta^3(\tilde s,
 t,\tilde x,y)$.\\
\textit{\underline{Term $R^3_\Delta(\tilde s,t,\tilde x,y) $.}} From
Lemma \ref{CTR_RESTES} applied to the process with frozen
coefficients,
one gets
\begin{align}
\label{CTR_R3}
R^3_\Delta(\tilde s,t,\tilde x,y)\le C(y_\Delta 
+\Delta^{1/4}).
\end{align}
\textit{\underline{Term $R_\Delta^2(\tilde s,t,\tilde x) $.}} Let us explain the leading ideas of the estimates below. Usually, it is easy to prove inequalities like $|Y_t-F(t,X^\Delta_t)|_{L_2}=O(\Delta^{1/2})$ (for a fixed $t$), but this not enough to control $R^2_\Delta$. To achieve our goal, we take advantage of the fact that the time $t$ is the stopping time $\tau^{\Delta,Y}$ which is likely close to $\tilde s$. Thus, $Y_{\tau^{\Delta,Y}}
 -F(\tau^{\Delta,Y},X_{\tau^{\Delta,Y}}^{\Delta})$ should be much smaller that $\Delta^{1/2}$ in $L_2$-norm.\\ 
Introduce for $0<\beta<\alpha<1/2,\ \tau_{\Delta^\beta}:=\inf \{s> \tilde s: |X_s^\Delta-\tilde x|\ge \Delta^\beta \}\wedge (\tilde s
+\Delta^\delta),\ \delta:=2\beta+\gamma,\ \gamma>0
  $. Clearly, one has
\begin{align*}
|R_\Delta^2(\tilde s,t,\tilde x)| \le& 2  \P_{\tilde s,\tilde
  x}\big[\tau^{\Delta,Y}\le t,\tau^\Delta=\tau^{\Delta,Y}, \tau^\Delta<\tau_{\Delta^\beta}
,
\\
&\qquad|Y_{\tau^{\Delta,Y}}
 -F(\tau^{\Delta,Y},X_{\tau^{\Delta,Y}}^{\Delta})|> y_\Delta\sqrt
 \Delta\big]
+2\P_{\tilde s,\tilde x}[\tau^\Delta\ge \tau_{\Delta^{\beta}},
\tau^\Delta\le t]\\
  :=& (R_\Delta^{21}+R_\Delta^{22})(\tilde s, t,\tilde x 
).
\end{align*}
Let us first deal with $R_\Delta^{21}(\tilde s, t,\tilde x)$.
By the Markov inequality, one has
\begin{align}
& R_\Delta^{21}(\tilde s,t,\tilde x)\le 2\Delta^{-1}y_\Delta^{-2}\E_{\tilde s,\tilde
x}\big[\I_{\tau^\Delta<\tau_{\Delta^\beta},\tau^{\Delta,Y}\le t,
\tau^\Delta=\tau^{\Delta,Y}
}|Y_{\tau^{\Delta,Y}}-F(\tau^{\Delta,Y},X_{\tau^{\Delta,Y}}^{\Delta})|^2\big]
. \label{First_Approx}
\end{align}
Note that since $\cD$ is of class $\H_2$, $F$ has the same
regularity, i.e. it is uniformly Lipschitz continuous in time,
its first space derivatives are uniformly Lipschitz continuous in space and $1/2$-Hölder continuous in time. Thus,
assuming up to a regularization procedure that $F\in
C^{1,2}([0,T]\times \R^d) $, Itô's formula yields for all $t\ge \tilde s$,
\begin{align}
F(t,X_{t}^\Delta)=&F(\tilde s,\tilde x)+\bint{\tilde
  s}^{t}\nabla F(s, X_s^\Delta) dX_s^\Delta\nonumber\\
  &+\bint{\tilde s}^{t}\bigl(\partial_s
F(s,X_s^\Delta)+\frac 12
\Tr(H_F(s,X_s^\Delta)\sigma\sigma^*(\phi(s),X_{\phi(s)}^\Delta))\bigr)ds\nonumber\\
:=& F(\tilde s,\tilde x)+\bint{\tilde s}^{t}\nabla F(
  s,  X_s^\Delta) \sigma(\phi(s),X_{\phi(s)}^\Delta) dW_s+R_F^\Delta(\tilde
  s,t,\tilde x)
  \label{NOT_POUR_PLUS_TARD}\\
=&Y_t+R_F^\Delta(\tilde s, t,\tilde x)+\int_{\tilde s}^{t}\bigl( \nabla F(
  s,  X_s^\Delta)\sigma(\phi(s),X_{\phi(s)}^\Delta)-[\nabla F\sigma](\tilde s,\tilde x)  \bigr) dW_s . \nonumber
\end{align}
From \A{A$_\theta^{'} $-1} and the assumptions on  $\cD$ one derives
$|R_F^\Delta|(\tilde s,t,\tilde x)\le C (t-\tilde s)$.
Thus, for any given stopping time $U\in[\tilde s, \tau_{\Delta^\beta}]
$, the working assumptions (i.e. smoothness of $\sigma,  F$) and standard
computations yield
\begin{align*}
\E[|F(U,X_U^\Delta)-Y_U|^2]\le C(\Delta^{2\beta+\delta}
+\Delta^{\delta(1+\theta)}).
\end{align*}
From \eqref{First_Approx} and the above control with $U=\tau^{\Delta,Y}\wedge \tau_{\Delta^\beta} $, one obtains
\begin{align}
\label{ctr_R21} R_\Delta^{21}(\tilde s,t,\tilde x)\le C
y_\Delta^{-2}\Delta ^{-1}(\Delta^{2\beta+\delta}
+\Delta^{\delta(1+\theta)}).
\end{align}
Let us now control $R_\Delta^{22}(\tilde s,t,\tilde x ) $. From
Lemmas \ref{contrprob} and \ref{CTR_Hitting}, for any $\eta>0$ we
write   
\begin{align}
R_\Delta^{22}(\tilde s, t,\tilde x)&\le \P_{\tilde s,\tilde
x}[\thb<\tilde
  s+\Delta^\delta]+\P_{\tilde s,\tilde x}[\tau^\Delta\land t\ge \tilde s+\Delta^\delta] \nonumber\\
&\le C_\eta\bigl(
\exp\left(-c\Delta^{2\beta-\delta}\right)+\Delta^{\alpha-\eta-\delta/2}+\Delta^{\delta/2}\bigr).
\label{ctr_R22}
\end{align}
Take now $\alpha=\frac{1+\frac\theta 2}{2(1+\theta)}<1/2,\
\eta=\frac{\theta}{16(\theta+1)},\ \gamma=\frac{1}{8(1+\theta)},\
y_\Delta=\Delta^{\theta/16}$. 
Check that for $\delta=2\beta+\gamma=2\alpha-4\eta$, one has $\delta=\frac{1+\theta/4}{1+\theta},\ \beta=\frac{7/8+\theta/4}{2(1+\theta)}<\alpha,\ 3\eta<\alpha$. Thus, $R_\Delta^{22}(\tilde s, t,\tilde x)=O(\Delta^\eta) $. In addition,
$y_\Delta^{-2}\Delta^{\delta(1+\theta)-1}=\Delta^{\theta/8},
y_\Delta^{-2}\Delta^{2\beta+\delta-1}=O(\Delta^{1/(8(1+\theta))}) $.
Hence, from \eqref{ctr_R21} and \eqref{ctr_R22}
\begin{align}
\label{CTR_R2} R_\Delta^{2}(\tilde s, t,\tilde x)&\le C\bigl(
\Delta^{1/(8(1+\theta))} +\Delta^{\theta/8}+
\Delta^{\theta/(16(\theta+1))}\bigr)\le C\Delta^{\theta/32}.
\end{align}
\textit{\underline{Term $R_\Delta^1(\tilde s,t,\tilde x) $.}} We
give an upper bound for $R_\Delta^{11}(\tilde s,t,\tilde x) $. The
term $R_\Delta^{12}(\tilde s,t,\tilde x) $ can be handled in the
same way. From the previous control on $R_\Delta^{22}(\tilde
s,t,\tilde x) $ and for the previous 
parameters, one gets 
\begin{align*}
R_\Delta^{11}(\tilde s,t,\tilde x)=&\P_{\tilde s,\tilde
x}[\tau^\Delta\le t,\tau^\Delta\neq
  \tau^{\Delta,Y}, \tau^\Delta <\tau_{\Delta^\beta}
  ]+ O(\Delta^\eta)
  \nonumber \\
=&\P_{\tilde s,\tilde x}[\tau^\Delta\le t,\tau^\Delta>
  \tau^{\Delta,Y}, \tau^\Delta <\tau_{\Delta^\beta}
  ]\nonumber\\
 & +\P_{\tilde s,\tilde x}[\tau^\Delta\le t,\tau^\Delta<
  \tau^{\Delta,Y}, \tau^\Delta<\tau_{\Delta^\beta}   ]
  + O(\Delta^\eta) 
  .
\end{align*}
Then, splitting the first probability according to $\Delta^{-1/2}(Y_{\tau^{\Delta,Y}})^-\le y_\Delta$ or not, and the second one according to $\Delta^{-1/2}F^-(\tau^\Delta,X_{\tau^{\Delta}}^{\Delta})\le y_\Delta$ or not, we obtain
\begin{align*}
&R_{\Delta}^{11}(\tilde s,t,\tilde x)\\
\le& \big(\P_{\tilde s,\tilde x}[\tau^{\Delta,Y}\le
  t,\Delta^{-1/2}(Y_{\tau^{\Delta,Y}})^-\le y_\Delta]\\
&+
\P_{\tilde s,\tilde  x}[\tau^\Delta\le
  t,\tau^\Delta>\tau^{\Delta,Y},  \tau^\Delta<\tau_{\Delta^\beta}
  ,
  \Delta^{-1/2}|Y_{\tau^{\Delta,Y}}-F(\tau^{\Delta,Y},X_{\tau^{\Delta,Y}}^{\Delta} ) |\ge y_\Delta]\big)\\
+& \big(\P_{\tilde s,\tilde
  x}[\tau^\Delta\le
  t,\tau^\Delta<\tau^{\Delta,Y}, \tau^\Delta<\tau_{\Delta^\beta},
\Delta^{-1/2}|Y_{\tau^{\Delta}}-F(\tau^{\Delta},X_{\tau^{\Delta}}^{\Delta}
) |\ge y_\Delta]\\
& +\P_{\tilde s,\tilde x}[\tau^{\Delta}\le
  t  ,\Delta^{-1/2}F^-(\tau^\Delta,X_{\tau^{\Delta}}^{\Delta})\le y_\Delta]\big)+C\Delta^\eta,
\end{align*}
for the previous function $(y_\Delta)_{\Delta> 0} $. Since we could obtain the same type of bound for
$R_\Delta^{12}(\tilde s,t,\tilde x) $, from Lemma \ref{CTR_RESTES} and
following the computations that gave \eqref{ctr_R21} we derive for the previous set of parameters
\begin{align}
\label{CTR_R1}
 R_\Delta^{1}(\tilde s,t,\tilde x)&\le
C
(y_\Delta^{-2}\Delta^{-1}(\Delta^{2\beta+\delta}+\Delta^{\delta(1+\theta)})
+\Delta^\eta+ y_\Delta+\Delta^{1/4})
\le C\Delta^{\theta/32} .
\end{align}
From \eqref{CTR_R1}, \eqref{CTR_R2}, \eqref{CTR_R3} we finally
obtain
$R_\Delta(\tilde s,t,\tilde x,y)=O(\Delta^{\theta/32})=o(1)$. The rest is uniform w.r.t. $(\tilde s,\tilde x,y)\in {\cal A}^{\alpha,\varepsilon}\times
\R^+ $.\\
\textbf{Step 4. Final step.} Plug the previous results in
\eqref{CTR_PSI}. We derive from \eqref{first_ctr}
\begin{align*}
\varPsi_\Delta(t,x,y)=\E_{x}&[\I_{
 A(t,\alpha,\varepsilon) }
  \varphi(
  X_\tha^\Delta
  )\\
&\times\P_{\tha,X_{\tha}^\Delta}[\tau^{\Delta,Y}\le
  t]
(1-H(y/
|\nabla F \sigma(\tha,
X_\tha^\Delta )
|) )
]
+o(1).
\end{align*}
Moreover, note that taking $y=0$ in the previous controls gives immediately
$$\P_{\tilde s,\tilde x}(\tau^{\Delta,Y}\le t)-\P_{\tilde s,\tilde x}(\tau^{\Delta}\le t)=o(1)$$ uniformly in $(\tilde s,\tilde x)\in{\cal A}^{\alpha,\varepsilon}$. Thus, we finally obtain
\begin{align*}
\varPsi_\Delta(t,x,y)=\E_{x}&[\I_{
 A(t,\alpha,\varepsilon) }
  \varphi(  X_\tha^\Delta  )\I_{\tau^{\Delta}\le
  t}(1-H(y/|\nabla F \sigma(\tha,X_\tha^\Delta )|))]
+o(1).
\end{align*}
Under continuity arguments as in step 1 (localization), we eventually get
\begin{align*}
\varPsi_\Delta(t,x,y)&=\E_{x}[\I_{\tau^{\Delta}\le
  t}\varphi( X_{\tau^\Delta}^\Delta) (1-H(y/
|\nabla F \sigma (\tau^\Delta, X_{\tau^\Delta}^\Delta)
|)) ]+o(1).
\end{align*}
We complete the proof using Lemma \ref{conv:proba:tau:xtau}:
\begin{align*}
  \varPsi_\Delta(t,x,y)\underset{\Delta\rightarrow 0}\rightarrow \E_{x}[\I_{\tau\le
  t}\varphi(X_\tau)
  (1-H(y/
|\nabla F \sigma(\tau,X_{\tau})  
|) )].
\end{align*}
\qed

\subsection{Proof of Lemmas \ref{lemme_siegmund}, \ref{CTR_Hitting} and \ref{CTR_RESTES}}
\hspace*{10pt}\\
\label{le_inter}

\textit{Proof of Lemma \ref{lemme_siegmund}.} We shall insist on the dependence of the exit times with respect to $x$, by setting $\tau^\Delta:=\inf\{t_i=i\Delta>0: W_{t_i}\ge x \}:=\tau^{\Delta}_x $ and analogously for $\tau=\tau_x$. Our proof relies on the following convergence (see equation (19) in Siegmund \cite{sieg:79}): if we set (for any $y,z\ge 0$)
$$D(z,y)=\P_0[W_{\tau_z^\Delta}-z \le y\sqrt{\Delta}]-H(y),$$
then
\begin{equation*}
  \lim_{z\Delta^{-1/2} \rightarrow +\infty} |D(z,y)|=0.
\end{equation*}
Using the monotonicity and the uniform continuity of $H(y)$, Dini's Theorem yields that the above limit is actually uniform with respect to $y\ge 0$. It follows
\begin{equation}
  \label{eq:sieg1}
 \sup_{y\ge 0,\ z\in[\Delta^{1/2-\varepsilon/3},\infty)}|D(z,y) |\underset{\Delta \rightarrow 0}{\rightarrow}0 .
\end{equation}
Additionnally, we have 
\begin{equation}
  \label{eq:sieg2}
 \sup_{x\ge 0,\ t\in[\Delta^{1-4 \varepsilon/3},T]}|\P_0(\tau_x^\Delta>t)-\P_0(\tau_x>t)|\underset{\Delta \rightarrow 0}{\rightarrow}0 .
\end{equation}
To prove this, we apply Theorem 3.4 in \cite{avik:07} which states that
\begin{equation*}
  \sup_{x \in \R}\E|\I_{M<x}-\I_{\hat M<x}|\le 3 (\sup_{m\in\R}f_M(m) \|M-\hat M\|_{L_p})^{\frac{p}{p+1}}
\end{equation*}
for any $p>0$ and for any random variables $M$ and $\hat M$, such that $M$ has a bounded density $f_M(.)$. Now, consider $M=\sup_{s\le t}W_s$ and $\hat M=\sup_{s=i\Delta\le t}W_s$. The density of $M$ is bounded by $2/\sqrt{2 \pi t}$. On the other hand, Lemma 6 in \cite{asmu:glyn:pitm:95} gives $\|M-\hat M\|_{L_p}\le C_p(T) \Delta^{1/2}$. Hence, we get for $t\ge \Delta^{1-4 \varepsilon/3}$, 
$$|\P_0(\tau_x^\Delta>t)-\P_0(\tau_x>t)|\le\E|\I_{\hat M<x}-\I_{M<x}| \le C_p(T) \Delta^{\frac{2 \varepsilon p}{3(p+1)}},$$
which leads to (\ref{eq:sieg2}).

We can now proceed to the proof of Lemma \ref{lemme_siegmund}, assuming that $x\ge \Delta^{1/2-\varepsilon}$. First, note that if $x/\sqrt t\ge \Delta^{-\varepsilon/3} \rightarrow +\infty$ as $\Delta\rightarrow 0$, $\P_0(\tau_x^\Delta\le t)$ and $\P_0(\tau_x\le t)$ are both $O_{pol}(\Delta)$. Thus, the difference in Lemma \ref{lemme_siegmund} converges to 0 as $\Delta\rightarrow 0$.

Suppose now that $x/\sqrt t\le \Delta^{-\varepsilon/3}$, hence $\sqrt t\ge x \Delta^{\varepsilon/3}\ge \Delta^{1/2-2\varepsilon/3}$, and write for $t\in\Delta \N^*$
\begin{align*}
P:=\P_0[\tau_x^\Delta>t,W_{\tau^\Delta_x}-x\le y\sqrt{\Delta}]
=\bint{0}^{+\infty}q_t^{x,\Delta}(0,x-z)\P_0[W_{\tau_z^\Delta}-z \le y\sqrt{\Delta}]dz
\end{align*}
where $q_{t}^{x,\Delta}(.,.)$ denotes the transition density of the Brownian motion discretely killed at level $x$. Introduce the partition $\R^+=[0,\Delta^{1/2-\varepsilon/3}) \cup[\Delta^{1/2-\varepsilon/3},+\infty)  $. Then,
\begin{align*}
  P&=R+\int_{\Delta^{1/2-\varepsilon/3}}^{+\infty} q_t^{x,\Delta}(0,x-z)D(z,y) dz +\P_0[\tau_x^\Delta>t]H(y)
\end{align*}
where $|R |\le  2\P_0[W_t\in[x-\Delta^{1/2-\varepsilon/3},x ]]\le \frac{2}{\sqrt{2\pi t}}\Delta^{1/2-\varepsilon/3}\le \frac{2}{\sqrt{2\pi}}\Delta^{\varepsilon/3}$ since $\sqrt t\ge \Delta^{1/2-2\varepsilon/3}$. Finally, taking advantage of the estimates \eqref{eq:sieg1} and \eqref{eq:sieg2} readily completes our proof.\qed


\textit{Proof of Lemma \ref{CTR_Hitting}.} We take $s=0$ for
notational simplicity. Introduce $\thb :=\inf\{t\ge 0:
X_t^\Delta \notin V_{\partial D_t}(\Delta^\beta)\}$ and for $\gamma>0 $
write from Lemma \ref{contrprob} and the notation of
\eqref{NOT_POUR_PLUS_TARD} (up to the same regularization procedure concerning $F$)
\begin{align*}
\P_x[\tau^\Delta \land T\ge \Delta^{2\beta}]=& \P_x[\inf_{0\le i\le\Delta^{2\beta -1}}  \big(F(0,x)+\int_{0}^{t_i}
  \nabla F(s, X_s^\Delta)\sigma(\phi(s),X_{\phi(s)}^\Delta)dW_s\\
&\quad+R_F^\Delta(0,t_i,x) \big) \ge 0,\thb\ge \Delta^{2\beta
+\gamma}]+O_{pol}(\Delta):=Q ,
\end{align*}
where under the assumptions of the Lemma, $|R_F^\Delta(0,t_i,x)|\le
C t_i  $ and $F(0,x)\le \Delta^\alpha$. For a given $r>0$, consider the event  ${\cal A}_r=\{\exists s\le T: |X^\Delta_s - X^\Delta_{\phi(s)}|\geq r\}$ where the increments of $X^\Delta$ between two close times are large: by Lemma \ref{contrprob}, it has an exponentially small probability. Hence, if we set
\begin{align*}
M_{u}:=& \int_{0}^{u}\nabla F(s, X_s^\Delta) \sigma
(\phi(s),X_{\phi(s)}^\Delta)d
  W_s:=B_{\langle M\rangle_u},\
  \tilde t_i=\langle M\rangle_{t_i},
\end{align*}
$B $ is a standard Brownian motion (on
a possibly enlarged probability space) owing to the Dambis, Dubbins-Schwarz Theorem, cf. Theorem
V.1.7 in \cite{revu:yor:99}. In addition, the above time
change is strictly increasing on the set $\cal A_r^c$ and $\langle M\rangle_t-\langle M\rangle_s\ge (t-s)a_0/2$ ($t\ge s$) up to taking $r$ small enough, because \A{A$_{\mathbf{\alpha}}^{'} $-2} is in force. It readily follows that
\begin{align*}
Q\le 
& \P_x[\inf_{0\le i\le\Delta^{2\beta+\gamma -1}}  (M_{t_i}+Ct_i)  \ge
    -\Delta^\alpha,\thb\ge \Delta^{2\beta +\gamma} ]+O_{pol}(\Delta)\\
\le & \P_x[\inf_{0\le i\le\Delta^{2\beta+\gamma -1}} (B_{\tilde t_i}+2Ca_0^{-1}\tilde t_i)\ge
    -\Delta^\alpha,\thb\ge \Delta^{2\beta +\gamma} ,\cal A_r^c]+O_{pol}(\Delta)\\
\le &\P_{x}[\inf_{0\le i\le\Delta^{2\beta+\gamma -1}}  (B_{\tilde t_i}+2Ca_0^{-1}\tilde t_i)\ge
  -\Delta^{\alpha},\thb  \ge \Delta^{2\beta+\gamma},\\
&\qquad\inf_{0\le s\le \langle M\rangle_{\Delta^{2\beta
  +\gamma}}} (B_s
+2Ca_0^{-1}s)
 \le
  -\Delta^{\alpha-\zeta},\cal A_r^c]+O_{pol}(\Delta)\\
&+\P_{x}[\thb \ge
  \Delta^{2\beta+\gamma},\inf_{0\le s\le \langle M\rangle_{\Delta^{2\beta
  +\gamma}} } (B_s+2Ca_0^{-1}s)\ge
  -\Delta^{\alpha-\zeta},\cal A_r^c],
\end{align*}
for $\zeta>0 $. Thus, from Lemma \ref{contrprob} and standard controls
\begin{align*}
Q&\le \P_{x}[\exists i: 0\le i \le \Delta^{2\beta+\gamma-1},\sup_{s\in[\tilde t_i,\tilde
      t_{i+1}]}|B_s-B_{\tilde t_i}+2Ca_0^{-1}(s-\tilde t_i) |\ge
      \Delta^{\alpha-\zeta}-\Delta^\alpha, \nonumber \\
& \thb\ge \Delta^{2\beta+\gamma}]+\P_{x}[\inf_{0\le s \le a_0
      \Delta^{2\beta+\gamma}/2}B_s\ge
      -\Delta^{\alpha-\zeta}-C 
\Delta^{2\beta
      +\gamma}]+O_{pol}(\Delta)\nonumber\\
& \le O_{pol}(\Delta)+C(\Delta^{\alpha-\zeta-\beta-\gamma/2}+\Delta^{\beta+\gamma/2 }). 
\end{align*}
Choose now $\gamma, \zeta $ s.t. $(\zeta+\frac\gamma 2)=\eta>0 $.
The proof is complete. \qed 

\textit{Proof of Lemma \ref{CTR_RESTES}.} 
Taking also $s=0$ for notational convenience, we write
\begin{align}
P:=&\P_x[\tau^\Delta\le
  t,\Delta^{-1/2}F^-(\tau^\Delta,X_{\tau^\Delta}^\Delta)\in[a,b]  ]\le  O_{pol}(\Delta)\nonumber \\
&+\bsum{i=1}^{ \lfloor
t/\Delta\rfloor}\E_{x}[\I_{\tau^\Delta>t_{i-1},X_{t_{i-1}}^\Delta
\in V_{\partial D_{t_{i-1}}}(r_0)}
\P_{\F_{t_{i-1}}}[\Delta^{-1/2}F^-(t_i,X_{t_i}^{\Delta})\in
  [a,b]
] ]
 \label{CTR_R}
\end{align}
using Lemma \ref{contrprob} for the last identity.\\
A Taylor formula gives:
$F(t_i,X_{t_i}^\Delta)
=F(t_{i-1},X_{t_{i-1}}^\Delta)+
\Sigma_{t_{i-1}}(W_{t_i}-W_{t_{i-1}})
+R_{t_{i-1},t_i}^\Delta:=\No_{t_{i-1}}+R_{t_{i-1},t_i}^\Delta
$ where $\Sigma_{t_{i-1}}= \nabla F \sigma(t_{i-1},X_{t_{i-1}}^\Delta),\ \E_{\F_{t_{i-1}}}[|R_{t_{i-1},t_i}^\Delta |^2]$ 
$\le C\Delta^2 $. Conditionally to $\F_{t_{i-1}}$, $\No_{t_{i-1}}$ has a Gaussian distribution \break $\No(F(t_{i-1},X_{t_{i-1}}^\Delta),\|\Sigma_{t_{i-1}} \|^2 \Delta).$\\
In addition, on the event $X_{t_{i-1}}^\Delta\in V_{\partial D_{t_{i-1}}(r_0)}, \
\|\Sigma_{t_{i-1}} \|^2\Delta \ge a_0\Delta $ and we obtain 
\begin{align*}
Q_{i-1}:=&\P_{\F_{t_{i-1}}}[F^-(t_i,X_{t_i}^\Delta) \in [a
    \Delta^{1/2},b \Delta^{1/2}]]\\
=&\P_{\F_{t_{i-1}}}[(\No_{t_{i-1}}+R_{t_{i-1},t_i}^\Delta)^-\in
  [a\Delta^{1/2},b \Delta^{1/2}]
  ]\\
   \le&  \P_{\F_{t_{i-1}}}[\No_{t_{i-1}}\in[-b\Delta^{1/2}-\Delta^{3/4},-a\Delta^{1/2}+\Delta^{3/4}]]\\
&\hspace{15pt}+\P_{\F_{t_{i-1}}}[|R_{t_{i-1},t_i}^\Delta |\ge  \Delta^{3/4},{X_{t_i}^\Delta\notin\cD_{t_i}}]\\
\le&
\P_{\F_{t_{i-1}}}[\No_{t_{i-1}}\in[-\Delta^{1/2}(b+\Delta^{1/4}),
-\Delta^{1/2}(a-\Delta^{1/4})]]\\
&\hspace{15pt}+C\Delta^{1/4}\exp\left(-c\bfrac{d(X_{t_{i-1}}^\Delta,\partial D_{t_{i-1}})^2 }{\Delta} \right)
\end{align*}
using the Cauchy-Schwarz inequality and Lemma \ref{contrprob} for the last inequality.
Hence, we derive from \eqref{CTR_R}
\begin{align*}
P\le  &\bsum{i=1}^{ \lfloor t/\Delta
  \rfloor}\E_{x}[\I_{\tau^\Delta>t_{i-1},X_{t_{i-1}}^\Delta\in V_{\partial D_{t_{i-1}}}(r_0)} \left(C\Delta^{1/4}\exp\left(-c\bfrac{d(X_{t_{i-1}}^\Delta,\partial D_{t_{i-1}})^2 }{\Delta} \right)\right. \\
+&\left.\bint{-\Delta^{1/2}(b+\Delta^{1/4})}^{-\Delta^{1/2}(a-\Delta^{1/4})}\exp\left(-\frac{(y-F(t_{i-1},X_{t_{i-1}}^{\Delta}))^2}{2\|\Sigma_{t_{i-1}}
  \|^2\Delta} \right)
\frac{dy}{(2\pi \Delta
  )^{1/2}\|\Sigma_{t_{i-1}} \|} \right)]+O_{pol}(\Delta).\end{align*}
We now upper bound the above integral on the event $\{\tau^\Delta>t_{i-1}
\}\subset\{F(t_{i-1},X_{t_{i-1}}^{\Delta})>0\}$. 
\begin{itemize}
\item If $y\le 0$, clearly one has $(y-F(t_{i-1},X_{t_{i-1}}^{\Delta}))^2\ge F^2(t_{i-1},X_{t_{i-1}}^{\Delta})$.
\item If $y\in(0,[\Delta^{1/2}(\Delta^{1/4}-a)]_+)$, one has $(y-F(t_{i-1},X_{t_{i-1}}^{\Delta}))^2\ge \frac 1 2 F^2(t_{i-1},X_{t_{i-1}}^{\Delta}) - y^2\ge \frac 1 2 F^2(t_{i-1},X_{t_{i-1}}^{\Delta}) - \Delta^{3/2}$.
\end{itemize}
Thus, we obtain that $P$ is bounded by
\begin{align*}C (b-a+\Delta^{1/4}) \bsum{i=1}^{ \lfloor t/\Delta
  \rfloor}\E_{x}[\I_{\tau^\Delta>t_{i-1},X_{t_{i-1}}^\Delta\in V_{\partial D_{t_{i-1}}}(r_0)} \exp(-c\bfrac{F^2(t_{i-1},X_{t_{i-1}}^\Delta)}{\Delta})]+O_{pol}(\Delta).
\end{align*}
The end of the proof is now achieved by standard computations done in \cite{gobe:meno:spa:04} p. 212 to 217. We only mention the main steps and refer for the details to the above reference. First, we replace the discrete sum on $i$ by a continuous integral, then we apply the occupation time formula to the distance process $(F(s,X_s^\Delta))_{s\leq \tau^\Delta} $ using the non characteristic boundary condition, as in the proof of Theorem \ref{FIRST_EXP}:
\begin{align*}
P\le  &C \frac{(b-a+\Delta^{1/4})}{\Delta}\int_0^t 
\E_{x}[\I_{\tau^\Delta>s,X_{s}^\Delta\in V_{\partial D_{s}}(r_0)} \exp(-c\bfrac{F^2(s,X_s^\Delta)}{\Delta})] ds +O_{pol}(\Delta)\\
\le &C \frac{(b-a+\Delta^{1/4})}{\Delta}\int_{-r_0}^{r_0}  \exp(-c\bfrac{y^2}{\Delta}) \E_{x}[L_{t\land \tau^\Delta}^y (F(.,X^\Delta_.))] dy+O_{pol}(\Delta).
\end{align*}
Then, we use \eqref{eq:tps:local}
to obtain $P\le C((b-a+\Delta^{1/4})$ which is our claim.
\qed

\begin{rem} Finally, we mention that if $\sigma\sigma^*$ is uniformly elliptic, the rest $R_{t_{i-1},t_i}^\Delta $ can be avoided and the result can be stated without the contribution $\Delta^{1/4}$. Indeed, we can directly exploit that the Euler scheme has conditionally a non degenerate Gaussian distribution and usual changes of chart associated to a parametrization of the boundary (see e.g. \cite{gobe:00}) give the expected result.
\end{rem}

\mysection{Extension to the stationary case}
\label{stat}
\subsection{Framework}
In this section we assume that the coefficients in \eqref{diff} are time
independent and that the mappings $b,\sigma $ are uniformly Lipschitz
continuous, i.e. $(X_t)_{t\ge 0}$ is the unique strong solution of
\begin{align*}
X_t=x+\bint{0}^{t}b(X_s)ds+\bint{0}^{t}\sigma(X_s)dW_s, t\ge 0,\ x\in\R^d.
\end{align*}
For a bounded domain $D\subset \R^d$, and given functions $f,g,k:\bar
D\rightarrow \R $, we are interested in estimating
\begin{align}
\label{rep:stat}
u(x):=\E_x[g(X_{\tau})Z_\tau+\bint{0}^{\tau} f(X_s)Z_s ds],\
Z_s=\exp(-\bint{0}^sk(X_r)dr),
\end{align}
where $\tau:=\inf\{t> 0:X_t\notin D \}$.

Adapting freely the previous notations for Hölder spaces to the elliptic setting, introduce for $\theta \in]0,1]$:
\begin{enumerate}
\item[\A{A$_\theta$}] 
   \begin{enumerate}
   \item[{\bf 1.}] \textit{Smoothness of the coefficients}. $b,\sigma \in \H_{1+\theta}$.
   \item[{\bf 2.}] \textit{Uniform ellipticity}. For some $a_0>0,\ \forall
   (x,\xi)\in\R^d\times \R^d,\ \xi^*\sigma \sigma^*(x)\xi$ $\ge a_0|\xi|^2 $.
  \end{enumerate}
\item[\A{D}\ ]\textit{Smoothness of the domain}. The bounded domain $D$ is
   of class $\H_{2}$.
\item[\A{C$_\theta$}]\textit{Other coefficients}. The boundary data $g\in
  \H_{1+\theta}$, $f,k\in \H_{1+\theta}$
 and $k\ge 0$.
\end{enumerate}
Note that under \A{A$_\theta$} and since $D$ is bounded, Lemma 3.1 Chapter III of \cite{frei:85}
yields $\sup_{x\in\bar D}\E_x[\tau]< \infty$. Thus, \eqref{rep:stat} is well defined under our current assumptions.

From Theorem 6.13, the final notes of Chapter 6 in \cite{gilb:trud:98} and Theorem 2.1 Chapter II in Freidlin
\cite{frei:85}, the Feynman-Kac representation in our elliptic setting writes
\begin{prop}[Elliptic Feynman-Kac's formula and 
esti\-ma\-tes]
\label{ell_FK}
\hspace*{15pt}\\
As\-sume \A{A$_\theta$}, \A{D}, \A{C$_\theta $} are in
force. Then, there is a unique solution in $\H_{1+\theta}\cap {\cal C}^2(D)$ to
\begin{align}
\label{EDP:stat_EQ}
\left\{
\begin{array}{c}
Lu-ku+f=0,\ {\rm in}\ D,\\
u|_{\deD}=g,
\end{array}
\right.
\end{align}
(where $L$ stands for the infinitesimal generator of $X$) and the solution is given by \eqref{rep:stat}. 
\end{prop}
In the following we denote by $F(x) $ the signed spatial distance to
the boundary $ \partial D$. Under \A{D}, $D$ satisfies the exterior and interior
uniform sphere condition with radius $r_0>0$ and  $F\in
\H_{2}(V_{\deD}(r_0))$ where $V_{\deD}(r_0):=\{x\in \R^d: {\bf d}(x,\deD)\le r_0\}$. Also, $F$ can be extended to a
$\H_{2}$ function preserving the sign. For more details on the
distance function, we refer to Appendix 14.6 in
\cite{gilb:trud:98}.
\subsection{Tools and results}
Below, we keep the previous notations concerning the Euler scheme. We also use the symbol $C$ for nonnegative
constants that may depend on $D,b,\sigma,g,f,k $ but not on $\Delta $
or $x$. We reserve the notation $c$ for constants also independent
of $D,g,f,k$. 

We recall a known result from Gobet and Maire \cite{gobe:mair:05} (Theorem 4.2) which provides an uniform bound for the $p$-th moment of $\tau^\Delta$:
\begin{align}
\label{int_ts_Euler}
\forall p\ge 1,\ \limsup_{\Delta \rightarrow 0}\sup_{x\in \bar D} \E_x[(\tau^\Delta)^p]<\infty.
\end{align}
Let us now state the main results of Section \ref{main_results} in our current framework.
\begin{prop}[Tightness of the overshoot]
\label{tension_STAT}
Assume \A{${\mathbf  {A}}_\theta $-2}, and that $D$ is of class $\H_2 $. Then, for some $c>0$,
$$\sup_{\Delta>0}\E_x[\exp(c[\Delta^{-1/2} F^-(X_{\tau^\Delta}^\Delta )]^2)]<+\infty.$$
\end{prop}
From the proof of Theorem \ref{loi_limite} and the estimate (\ref{int_ts_Euler}) we derive:
\begin{thm}[Joint limit laws associated to the overshoot] \label{loi_limite_STAT}
Assume \A{${\mathbf {A}}_\theta$}, and that $D$ is of class
$\H_2$. Let $\varphi$ be a continuous function with compact
support. With the notation of Theorem \ref{loi_limite}, for all $x\in D,\ y\ge 0$,
\begin{align*}
&&\E_x[Z_{\tau^\Delta}^\Delta\varphi(
X_{\tau^\Delta}^\Delta
)\I_{F^-(X_{\tau^\Delta}^\Delta)\ge
y\sqrt
  \Delta}]\underset{\Delta \rightarrow 0}{\longrightarrow}
\E_x\bigl[Z_{\tau}\varphi(X_\tau)\bigl(1-H(y/| \nabla F\sigma 
  (X_\tau)|)\bigr)\bigr].
\end{align*}
\end{thm}

\subsection{Error expansion and boundary correction}
For notational convenience introduce for $x\in D$,
\begin{align*}
u(D)&=\E_x(g(X_{\tau})Z_{\tau}+\int_0^{\tau}Z_sf(X_s)ds),\\
u^\Delta(D)&=\E_x(g(
X^{\Delta}_{\tau^{\Delta}}
)Z^{\Delta}_{\tau^{\Delta}}+\int_{0}^{\tau^{\Delta}}
  Z^{\Delta}_{\phi(s)} f(X^{\Delta}_{\phi(s)}) ds).
\end{align*}
The second quantity is well defined owing to (\ref{int_ts_Euler}).
\begin{thm}[First order expansion]
\label{dev_STAT}
Under \A{A$_\theta$}, \A{D}, \A{C$_\theta$}, for $\Delta$
small enough and with the notation of Theorem \ref{dev}
\begin{align*}& \ellErr{}=u^\Delta(D)-u(D)\\
&=c_0 \sqrt \Delta \E_x(Z_\tau(\nabla u -\nabla g)(X_{\tau})\cdot \nabla F(X_{\tau}) | \nabla F \sigma (X_\tau)|)+o(\sqrt\Delta).
  \end{align*}
\end{thm}
Define now $D^\Delta=\{x\in D: {\bf d}(x,\deD)> c_0 \sqrt \Delta
| \nabla F \sigma (x)|\}$. Introduce
$\hat\tau^{\Delta}=\inf\{t_i> 0:X^\Delta_{t_i}\in D^\Delta\}$. Set
\begin{align*}
& u^\Delta(D^\Delta)= \E_x[g(X^{\Delta}_{\hat\tau^{\Delta}})Z^{\Delta}_{\hat\tau^{\Delta}}+\int_{0}^{\hat\tau^{\Delta}} Z^{\Delta}_{\phi(s)} f(X^{\Delta}_{\phi(s)}) ds].
\end{align*}
One has:
\begin{thm}[Boundary correction]
\label{corr_STAT}
Under \A{A$_\theta $}, \A{D}, \A{C$_\theta$} and assuming additionally $ \nabla F(.) | \nabla F \sigma (.)|$ is in ${\cal C}^{2}$, then 
for $\Delta$ small enough one has
$$u^\Delta(D^\Delta)-u(D)=o(\sqrt \Delta).$$
\end{thm}
\subsection{Proofs}
Note carefully that all the constants appearing in the error analysis for the parabolic
case have at most linear growth w.r.t the fixed final time $T$.
Estimate \eqref{int_ts_Euler} allows to control uniformly the integrability of these constants in our current framework.
Thus, since the arguments remain the same, we only give below sketches of the proofs.

\textit{Proof of Proposition \ref{tension_STAT}.} It is sufficient to
prove that there exist constants $\tilde c>0$ and $C$ s.t. $\forall A\ge 0,\
\sup_{\Delta >0}\P_x[F^{-}(X_{\tau^\Delta}^\Delta) \ge
    A\Delta^{1/2}]\le C\exp(-\tilde cA^2)$. Then any choice of $c<\tilde c$ is valid. For
  $x\in D$, we write
\begin{align*}
 P&:=\P_x[F^-(X_{\tau^\Delta}^\Delta)\ge A\Delta^{1/2}]\\
&=\bsum{i\in
    \N^*}^{}\E[\I_{\tau^\Delta>t_{i-1}}\I_{\tau_{t_{i-1}}^\Delta<t_i}\P[F^-(X_{t_i}^\Delta)\ge A\Delta^{1/2}|\F_{\tau_{t_{i-1}}^\Delta}]]
\end{align*}
where $\tau_{t_{i-1}}^\Delta:=\inf\{s\ge t_{i-1}: X_s^\Delta\notin D
\}$. From Lemma \ref{contrprob}, we get
\begin{align*}
P\le C \exp(-\tilde cA^2)\bsum{i\in\N^*}^{}\P[\tau^\Delta>t_{i-1},\tau_{t_{i-1}}^\Delta<t_i].
\end{align*}
Lemma 16 from \cite{gobe:meno:spa:04} remains valid under our current
assumptions and yields 
$$P\le
C\exp(-\tilde cA^2)\sum_{i\in\N^*}^{}\E[\I_{\tau^\Delta>{t_{i-1}}}(\P[X_{t_i}^\Delta
  \notin D|\F_{t_{i-1}}]+O_{pol}(\Delta))].$$
On the one hand, $\sum_{i\in\N^*}^{}\I_{\tau^\Delta>{t_{i-1}}}\I_{X_{t_i}^\Delta   \notin D}=\I_{\tau^\Delta<\infty}=1$ owing to \eqref{int_ts_Euler}. On the other hand, we have $\sum_{i\in\N^*}^{}\P_x[\tau^\Delta>t_{i-1}]=
\Delta^{-1}\E_x[\tau^\Delta]\le C/\Delta $ using \eqref{int_ts_Euler} again. Finally, we obtain
that $P\le C\exp(-\tilde cA^2) $ which concludes the proof.\qed

\textit{Proof of Theorem \ref{dev_STAT}.} Similarly to the proof of
Theorem \ref{FIRST_EXP} we suppose first that $u\in \H_{3+\theta}$. The general case can be deduced as in the parabolic case using suitable Schauder estimates, given in the final notes of Chapter 6 in \cite{gilb:trud:98}, see also our Appendix.

In this simplified setting, keeping the notations introduced in the proof of Theorem \ref{FIRST_EXP}, we obtain
\begin{align}
\ellErr{}\= Z_{\tau^\Delta}^\Delta(\nabla u- \nabla g)(\pi_{\partial D}(X_{\tau^\Delta}^\Delta))\nabla F(X_{\tau^\Delta}^\Delta)F^-(X_{\tau^\Delta}^\Delta)\nonumber\\
+ \bigl( \bsum{i\in\N}^{}\I_{t_i<\tau^\Delta}\bigl[\I_{A_{t_i}^\varepsilon} O(\Delta^{\frac{1+\theta}2})\label{R_PRES}\\
+\I_{(A_{t_i}^\varepsilon)^C}\I_{\forall s\in [t_i,t_{i+1}], X_s^\Delta\in B(X^\Delta_{t_i},\Delta^{\frac 12(1-\varepsilon)})}
 (u(X_{t_{i+1}}^\Delta)Z_{t_{i+1}}^\Delta-u(X_{t_i}^\Delta)Z_{t_i}^\Delta\nonumber\\+ Z_{t_i}^\Delta f(X_{t_i}^\Delta)\Delta) \bigr]\bigr)\I_{\tau^{r_0}>\tau^\Delta}. \label{R_LOIN}
\end{align}
Since the constant in \eqref{eq:tps:local} depends linearly on time, the contribution associated to the remainder \eqref{R_PRES} can be
bounded by
$C\Delta^{\frac{3+\theta-\varepsilon}2}$ $\times(\Delta^{-1}\E_x[\tau^{\Delta}])$. From \eqref{int_ts_Euler},
this quantity is a $O(\Delta^{\frac{1+\theta-\varepsilon}2})=o( \Delta^{\frac12}) $ for $\varepsilon $ small enough. Similarly to \eqref{DEP_DER_VARIEES} the term \eqref{R_LOIN} can be bounded by
\begin{align*}
&\E\bigl[\bigl(\bsum{i\in\N}^{}\I_{t_i<\tau^\Delta}\I_{(A_{t_i}^\varepsilon)^C
}
O(\Delta^2\{1+|u|_\infty+|\nabla u|_\infty+|D^2 u|_\infty+|D^3 u|_\infty\}\\
&\hspace{4.5cm}+\Delta^{\frac{3+\theta}2}[D^3 u]_{x,\theta}
)
\bigr) \bigr] +O_{pol}(\Delta)\\
&\le C\Delta^{\frac{1+\theta}2}\E[\tau^\Delta]=o(\Delta^{1/2}).
\end{align*}
 We eventually derive the result as in Section \ref{main_results}.
\qed

Theorem \ref{corr_STAT} can be proved as Theorem \ref{corr}, using a
sensitivity result analogous to Theorem 2.2 in \cite{cost:elka:gobe:06} for elliptic problems, see e.g. Simon \cite{simo:80}. We skip the details.

\mysection{Numerical results}
\label{resnum}
The numerical behavior of the correction of Theorem \ref{corr} had
already been illustrated for the killed case in Section 3 of
\cite{meno:06}. Additional tests are presented in \cite{gobe:09}. We now focus on the stopped case with the following example.  Take $d=3$ and introduce the following diffusion process
\begin{eqnarray}
dX_t&=&b(X_t)dt+\sigma(X_t)dW_t, \ \forall x \in \R^3, \ b(x)=\left(
x_2\ x_3\ x_1 
\right)^*,\nonumber\\
\sigma(x)&=&\left(\begin{array}{ccc}(1+|x_3|)^{1/2} &0 &0\\
\frac 12 (1+|x_1|)^{1/2}&\left(\frac34\right)^{1/2}(1+|x_1|)^{1/2}&0\\
0&\frac 12 (1+|x_2|)^{1/2}&\left(\frac34\right)^{1/2}(1+|x_2|)^{1/2}
\end{array} \right),\nonumber \\
\label{eq_ex}\
\end{eqnarray}
and $X_0$ to be specified later on. Set $D=B(0,2)$. We consider an elliptic problem. Starting from a given function $u(x)=x_1x_2x_3$ defined on $\bar D$, we derive the PDE of type \eqref{EDP:stat_EQ} associated to \eqref{eq_ex} satisfied by $u$ by taking $g=u|_{\partial D}$, setting  $f=-Lu$ where $L$ stands for the infinitesimal generator of  $X$ in \eqref{eq_ex} and $k=0$. 
One can easily check that $-f(x)=x_2^2x_3+x_3^2x_1+x_1^2x_2+\frac 12[x_3(1+|x_1|)^{1/2}(1+|x_3|)^{1/2} +x_1\left(\frac 34\right)^{1/2}(1+|x_1|)^{1/2}(1+|x_2|)^{1/2}] $.
Thus we have an explicit expression for the solution of \eqref{EDP:stat_EQ}.

For $x_0$ s.t. $(x_0^i)_{1\le i\le 3} \in \{-.7,-.3,.3,.7\} $, we take $N_{MC}=10^6$ sample paths for the Monte Carlo simulation and let $\Delta$ vary in $\{.01,.05,.1 \} $.
For all the computations, the size of the 95\% confidence interval always varies in $[1.5\times 10^{-3},  2\times 10^{-3}]$. For the absolute value of the absolute and relative errors over the $3\times 4^3=192$ points of the spatial grid, we report the results in Table \ref{table:1}. These results for the correction seem to indicate that the remainder $o(\Delta^{1/2}) $ in Theorem \ref{corr_STAT} is actually a $O(\Delta)$. This will concern further research. 

\begin{table}[h]
\begin{center}
\begin{tabular}{|c|c|c|}
\hline
 $\Delta$ & Without correction & In the corrected domain\\
\hline
.1 &    0.169 (199\%)  &        0.0220 (24.4\%)
                     \\
.05 &     0.114 (133\%) &       0.0115  (13.1\%)       			\\  
.01 &	0.0471 (54.7\%) &       0.0026  (2.98\%)  \\\hline
\end{tabular}
\end{center}
\caption{Supremum of the absolute error for the Euler scheme (relative error in $\%$ in parenthesis)\label{table:1}}
\end{table}

In Tables \ref{table:2} and \ref{table:3}, we also report the results obtained for the spatial points $x_0=(-.7 , .3, .7) $ and $x_0=(-.7,.7,-.7) $.
\begin{table}[h]
\begin{center}
\begin{tabular}{|c|c|c|}
\hline
 $\Delta$ & Without correction & In the corrected domain\\
\hline
.1 & -.0913+/- .0019 & -.1477 +/- .0016 \\
.05& -.1051 +/- .0018& -.1465+/- .0016   \\  
.01& -.1282 +/- .0017& -.1476+/- .0016  \\
\hline
\end{tabular}
\end{center}
\caption{Estimated value at $x_0=(-.7 , .3, .7) $ (with 95\% confidence interval). True value $u(x_0)=-.147 $. \label{table:2}}
\end{table}
\begin{table}[h]
\begin{center}
\begin{tabular}{|c|c|c|}
\hline
 $\Delta$ & Without correction & In the corrected domain\\
\hline
.1 &  .5368 +/- .0019& .3866 +/- .0016\\
.05&  .4648 +/- .0018 &  .3634 +/- .0016 \\  
.01&  .3851 +/- .0016 & .3473 +/- .0016 \\
\hline
\end{tabular}
\end{center}
\caption{Estimated value at $x_0=(-.7,.7,-.7) $ (with 95\% confidence interval). True value $u(x_0)=.343 $. \label{table:3}}
\end{table}

\begin{figure}[h]
\begin{center}
\includegraphics[scale=.52]{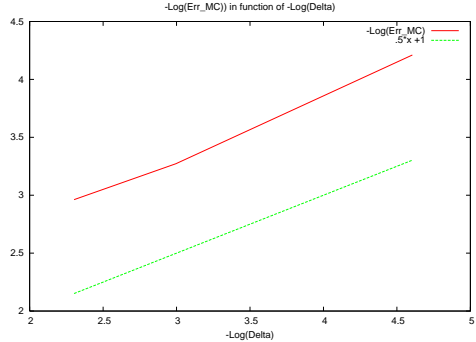}
\end{center}
\caption{Error for the Monte Carlo method (without correction) as a function of $\Delta$, in logarithmic scales. Evaluation at $x_0=(-.7, .3 , .7)$.\label{fig:4}}
\end{figure}

Eventually, for the Monte Carlo method, taking $x_0=(-.7, .3 , .7)$ and the previous values of $\Delta$, in Figure \ref{fig:4} we plot $-\log({\rm Err}_{MC}) $ in function of $-\log (\Delta)$, where ${\rm Err}_{MC}:=\left\{\frac{1}{MC}\bsum{i=1}^{MC}\left( g(X_{\tau^{\Delta,i}}^{\Delta,i})+\bint{0}^{\tau^{\Delta,i}}f(X_{\phi(s)}^{\Delta,i})ds\right)\right\}-u(x_0) $. The curve is quite close to a right line with slope $1/2$ as it should from Theorem \ref{dev_STAT}.

\mysection{Conclusion}
We have proposed and analysed a boundary correction procedure to simulate stopped/killed diffusion processes. This is valid for non-stationary and stationary problems, in time-dependent or time-independent domains. The resulting scheme is elementary to implement and its numerical accuracy is very good in our experiments. The proof relies on new asymptotic results regarding the renormalized overshoots.\\
To conclude, we note that the boundary correction procedure is very generic and could be at least formally extended to general Itô processes of the form $dX_t=b_t dt+\sigma_t dW_t$. In that case, the smaller domain would be defined $\omega$ by $\omega$ replacing $ \nabla  F(t,x)\sigma (t,x)$ by $\nabla F(t,X_t)\sigma_t$. Even if our current proof relies on Markovian properties, we conjecture that the correction should once again give a $o(\sqrt \Delta)$ independently of the Markovian structure. Numerical tests in \cite{gobe:09} support this conjecture, which will be addressed mathematically in further research.

\appendix

\mysection{Proof of Theorem {\ref{FIRST_EXP}} in the general setting}

In this section, we detail how the proof of Section \ref{main_results}
has to be modified under the assumptions of Theorem \ref{dev}, i.e. for $g\in\H_{1+\theta}$ and without compatibility condition so that $u\in \H_{1+\theta} $. Actually, $u$ is smooth inside the domain but high order derivatives may explode close to the boundary. These features have to be accurately quantified to show that the induced singularities are integrable. 

\subsection{Preliminary notation and controls}

Introduce the parabolic distance $\pd$: for $(s,x),(t,y)\in \bar{\cD},\ \pd((s,x),(t,y))=\max(|s-t|^{1/2},|x-y|)$. We also denote for a closed set $\cA \in \bar \cD$ and $(s,x)\in\cD,\ \pd((s,x),\cA) $ the parabolic distance of $(s,x) $ to $\cA $. Note that $\pd((s,x),\PD \cap\{v\ge s \})\ge \min(F(s,x),\sqrt{T-s})$, so that we obtain the easy inequality:
\begin{equation}
  \label{eq:pd}
  \frac1{\pd((s,x),\PD \cap\{v\ge s \})}\le \frac1{F(s,x)}+\frac1{\sqrt{T-s}}.
\end{equation}
Under our current assumptions, for some constant $C>0$, we have
\begin{align}
\label{CTR_SCHAUDER1}
 |D^2 u(s,x)|+|D^3 u(s,x) | \le C \pd((s,x),\PD
\cap\{v\ge s \})^{-2};\\
\text{ for } (t,y)\neq(s,x),\quad  \frac{|D^3 u(s,x)-D^3 u(t,y)|}{\pd((s,x),(t,y))^\theta}\nonumber\\
\le C [\pd((s,x),\PD \cap\{v\ge s \})\wedge \pd((t,y),\PD \cap\{v\ge
t \})]^{-2-\theta}; \label{CTR_SCHAUDER2}\\
\text{ for } t\neq s,\quad\frac{|D^2 u(s,x)-D^2 u(t,x)|}{|t-s|^{(1+\theta)/2}}\nonumber\\
\le C[\pd((s,x),\PD \cap\{v\ge s \})\wedge \pd((t,x),\PD \cap\{v\ge
t \})]^{-2-\theta}.\label{CTR_SCHAUDER3}
\end{align}
The above constant $C$ is uniform w.r.t. $(s,x)\in\D, (t,y)\in \D$ or $(t,x)\in \D$.
These inequalities are obtained with the interior Schauder estimates for the PDEs satisfied by the partial derivatives $(\partial_{x_i}u)_{1\le i\le d}$, see Theorem 4.9 in \cite{lieb:96}.


We first state an important proposition for the error analysis
with possibly explosive controls as in \eqref{CTR_SCHAUDER1}-(\ref{CTR_SCHAUDER2})-(\ref{CTR_SCHAUDER3}) for the derivatives.
Namely, under our current regularity assumptions, in order to perform a Taylor
expansion we have to work with interior points located in small balls, which distance to the boundary is uniformly bounded from below within the ball. The next proposition states that this is the case if the ball centers are "far enough" from the side of $\D$.
\begin{prop} \label{int_shift}
Assume $\cD\in \H_2 $ and take $\varepsilon\in]0,1[$. For all $(t,x)\in \bar \D\cap V_{\dcD}(r_0/2)\backslash V_{\dcD}(2\Delta^{1/2(1-\varepsilon)})$ ($r_0$ is defined in Section \ref{time domains}), one has for $ \forall y\in B(x,\Delta^{1/2(1-\varepsilon)})$ and $s\in[t,t+\Delta]$
\begin{align*}
F(s,y)\ge \frac 14 F(t,x)
\end{align*}
for $\Delta $ small enough (uniformly in $t,x,s,y$). In particular, $y$ belongs to $D_s$.
\end{prop}

\textit{Proof.} 
Since $F\in \H_2$, one has
\begin{align*}
F(s,y)\ge &F(t,x) -C\Delta+
 \langle \nabla F(t,x),y-x \rangle-C\Delta^{1-\varepsilon}
.
\end{align*}
The norm of $\nabla F(t,x)$ equals 1, since $\nabla F(t,x)$ is the unit inward normal vector at the closest point of $x$ on $\partial D_t$. Therefore, for $\Delta$ small enough and using $\frac 12 F(t,x)\ge \Delta^{\frac12(1 -\varepsilon)}$, we have $$F(s,y)\ge F(t,x) -\frac 32 \Delta^{\frac12(1 -\varepsilon)}\ge \frac 14 F(t,x),$$ which is the expected inequality.
\qed

We are now in a position to deduce useful local upper bounds for the derivatives of $u$ and their H\"older-norms, under the assumptions of Theorem \ref{FIRST_EXP}.
\begin{cor} \label{cor:schauder} Take $\varepsilon\in]0,1[$. There exists a constant $C>0$ such that for $\Delta$ small enough, for all $(t,x)\in \bar \D\backslash V_{\dcD}(2\Delta^{1/2(1-\varepsilon)})$, for all $(y,z)\in B(x,\Delta^{1/2(1-\varepsilon)})$ and $(r,s)\in[t,t+\Delta]$, we have
\begin{align}
\label{CTR_SCHAUDER1bis}
 |D^2 u(s,y)|+|D^3 u(s,y) | &\le \frac C{F^2(t,x)}+\frac C{T-t};\\
\text{ for } y\neq z,\quad  \frac{|D^3 u(s,y)-D^3 u(s,z)|}{|y-z|^\theta}&\le \frac C{F^{2+\theta}(t,x)}+\frac C{(T-t)^{1+\theta/2}};\label{CTR_SCHAUDER2bis}\\
\text{ for } r\neq s,\quad \frac{|D^2 u(r,y)-D^2 u(s,y)|}{|r-s|^{(1+\theta)/2}}&\le \frac C{F^{2+\theta}(t,x)}+\frac C{(T-t)^{1+\theta/2}}.\label{CTR_SCHAUDER3bis}
\end{align}
\end{cor}
\textit{Proof.} Note that if $(t,x)\in \bar \D\backslash V_{\dcD}(2\Delta^{1/2(1-\varepsilon)})$, we have $T-t\ge 4\Delta^{1-\varepsilon}$.\\
{\em Estimate (\ref{CTR_SCHAUDER1bis}).} In view of (\ref{CTR_SCHAUDER1}) and (\ref{eq:pd}), the upper bound of $|D^2 u(s,y)|+|D^3 u(s,y) |$ is equal to $\frac C{F^2(s,y)}+\frac C{T-s}$. On the one hand, by easy computations, we prove
$$\frac{1}{T-s}\le \frac{1}{T-t} \frac{T-t}{T-t-\Delta}\le \frac{1}{T-t} \frac{1}{1-\Delta^{\varepsilon}/4}\le \frac{C}{T-t}$$
for $\Delta$ small enough. On the other hand, we have 
$$\frac{1}{F(s,y)}\le \frac{C}{F(t,x)}.$$
Indeed, if $x$ is far from $D_t$ (and thus $y$ far from $D_s$), both terms $F(s,y)$ and $F(t,x)$ are bounded from above and from below. In the other case when $(t,x)\in \bar \D\cap V_{\dcD}(r_0/2)\backslash V_{\dcD}(2\Delta^{1/2(1-\varepsilon)})$, Proposition \ref{int_shift} yields $\frac{F(s,y)}{F(t,x)}\ge \frac 14$. Therefore, the upper bound (\ref{CTR_SCHAUDER1bis}) readily follows.\\
{\em Estimates (\ref{CTR_SCHAUDER2bis}) and (\ref{CTR_SCHAUDER3bis})}. They are proved following the same arguments, the details of which are left to the reader.
\qed
\subsection{Error analysis}
Recall from the previous proof of Theorem \ref{FIRST_EXP} that the main term to analyse is
\begin{align*}
e_{22}^\Delta&
\=\bigl(\sum_{0\le t_i< T}\1_{t_i<\tau^\Delta}\I_{(A_{t_i}^\varepsilon)^C}\I_{\forall s\in[t_i,t_{i+1}],\ X_s^\Delta \in B(X_{t_i}^\Delta, \Delta^{\frac12(1-\varepsilon)})}\left[
 u(t_{i+1},X^\Delta_{t_{i+1}})Z^\Delta_{t_{i+1}}\right.\nonumber \\
&\qquad \left.-u(t_{i},X^\Delta_{t_{i}})Z^\Delta_{t_i} + Z^\Delta_{t_i} f(t_i,X^\Delta_{t_{i}})\Delta\right]\bigr)\I_{\tro>\tau^\Delta\land T}\nonumber \\
&=\bigl(\sum_{0\le t_i< T-4\Delta^{1-\varepsilon}} \cdots\bigr)\I_{\tro>\tau^\Delta\land T}+\bigl(\sum_{T-4\Delta^{1-\varepsilon}\le t_i<T}\cdots\bigr)\I_{\tro>\tau^\Delta\land T}:= e_{221}^\Delta+e_{222}^\Delta,
\end{align*}
where we have just splitted the summation on $t_i$. 

\textit{Control of $e_{221}^\Delta $.} The idea is to perform a stochastic expansion of $ u(t_{i+1},X^\Delta_{t_{i+1}})Z^\Delta_{t_{i+1}} \break - u(t_{i},X^\Delta_{t_{i}})Z^\Delta_{t_i} + Z^\Delta_{t_i} f(t_i,X^\Delta_{t_{i}})\Delta$ as in \eqref{PREAL_DEP_DER_VARIEES}. Under our current assumptions, the difference comes from the high-order derivatives that are no more uniformly bounded or uniformly H\"older but only locally, with local estimates given in Corollary \ref{cor:schauder}. Thus, following the same computations that have led to \eqref{DEP_DER_VARIEES}, we obtain
\begin{align}
e_{221}^\Delta
\=& \bigl(\sum_{0\le t_i< T-4\Delta^{1-\varepsilon}}\1_{t_i<\tau^\Delta}\I_{(A_{t_i}^\varepsilon)^C}\I_{\forall s\in[t_i,t_{i+1}],\ X_s^\Delta \in B(X_{t_i}^\Delta, \Delta^{\frac12(1-\varepsilon)})}\bigl[\nonumber\\
&\quad O\bigl ((\Delta^2+\Delta |X_s^\Delta-X_{t_i}^\Delta|^2)(\frac 1{F^2(t_i,X_{t_i}^\Delta)}+\frac 1{T-t_i})\bigr)\nonumber\\
&+O\bigl ((\Delta^{1+\frac{1+\theta}{2}}+\Delta |X_s^\Delta-X_{t_i}^\Delta|^{1+\theta})(\frac 1{F^{2+\theta}(t_i,X_{t_i}^\Delta)}+\frac 1{(T-t_i)^{1+\theta/2}})\bigr)\nonumber\\
&+\bar M_{t_i,t_{i+1}}\bigr]\bigr)\I_{\tro>\tau^\Delta\land T}.\label{DEP_DER_VARIEES_MOD}
\end{align}
The derivatives appearing in $(\bar M_{t_i,t_{i+1}})_{0\le t_i<T}$ (see equations \eqref{THE_DETAILS} and \eqref{DEP_DER_VARIEES}) are controlled by (\ref{CTR_SCHAUDER1bis}) on $(A_{t_i}^\varepsilon)^C $. The control of \eqref{CTR_E22_0} remains valid for the $(\bar M_{t_i,t_{i+1}})_{0\le t_i<T} $ that yields a negligeable contribution. It follows that
\begin{align*}
|e_{221}^\Delta|\< C\Delta^{\frac{1+\theta}2}\bigl(\sum_{0\le t_i< T-4\Delta^{1-\varepsilon}}&\1_{t_i<\tau^\Delta}\I_{F(t_i,X_{t_i}^\Delta)\ge 2\Delta^{\frac12(1-\varepsilon)}}\Delta\bigl[\\
&\frac 1{F^{2+\theta}(t_i,X_{t_i}^\Delta)}+\frac 1{(T-t_i)^{1+\theta/2}}
\bigr]\bigr)\I_{\tro>\tau^\Delta\land T}.
\end{align*}
Standard computations show that 
$$\Delta^{\frac{1+\theta}2} \sum_{0\le t_i< T-4\Delta^{1-\varepsilon}}\frac \Delta{(T-t_i)^{1+\theta/2}}\le \Delta^{\frac{1+\theta}2}\int_0^{T-4\Delta^{1-\varepsilon}+\Delta} \frac{dt}{(T-t)^{1+\theta/2}}=O( \Delta^{\frac 12+\frac{\theta\varepsilon}2}),$$
which implies
$$|e_{221}^\Delta|\<C\Delta^{\frac{1+\theta}2}\biggl(\bint{0}^{T\wedge
\tau^\Delta}\I_{F(\phi(t),X_{\phi(t)}^\Delta)\in
[2\Delta^{1/2(1-\varepsilon)},r_0/2]
}F(\phi(t),X_{\phi(t)}^\Delta)^{-2-\theta}dt\biggr)+O( \Delta^{\frac 12+\frac{\theta\varepsilon}2}).$$
Adapting the previous analysis of Section \ref{main_results} for the term $e_{211}^\Delta $, we get
\begin{align*}
|e_{221}^\Delta|&\<C\Delta^{\frac{1+\theta}2}\biggl(\bint{0}^{T \wedge \tau^\Delta
}\I_{F(t,X_t^\Delta)\in [\Delta^{1/2(1-\varepsilon)},3 r_0/4]
}F(t,X_{t}^\Delta)^{-2-\theta}dt\biggr)+O( \Delta^{\frac 12+\frac{\theta\varepsilon}2})\\
&\<C\Delta^{\frac{1+\theta}2}\biggl(\bint{\Delta^{1/2(1-\varepsilon)}}^{3r_0/4}
y^{-2-\theta}L_{T\wedge
\tau^\Delta}^y(F(.,X_{.}^\Delta)) dy\biggr)+O( \Delta^{\frac 12+\frac{\theta\varepsilon}2}),
\end{align*}
using Lemma \ref{contrprob} for the last but one inequality, and the
occupation time formula for $F(t,X_t^\Delta) $ for the last one
(recall that $\sigma $ is uniformly elliptic).\\ Finally using \eqref{eq:tps:local}, one gets
\begin{align*}
|e_{221}^\Delta|&\le C\Delta^{\frac{1+\theta}2}\left(
\bint{\Delta^{1/2(1-\varepsilon)}}^{3r_0/4} y^{-2-\theta}
(y+\Delta^{1/2})dy\right)+O( \Delta^{\frac 12+\frac{\theta\varepsilon}2})
\le C\Delta^{\frac 12+\frac{\theta\varepsilon}2}=o(\Delta^{1/2}).
\end{align*}

\textit{Control of $e_{222}^\Delta$}. Apply a Taylor formula with integral rest at order one in space. The $\theta$-H\"older continuity in space of $\nabla u $ and the $(1+\theta)/2$-H\"older continuity in time of $u $ directly give a contribution in $O(\Delta^{1/2+\theta/2-\varepsilon})=o(\Delta^{1/2})$ for $\varepsilon $ small enough.
This completes the proof.\qed

\bibliography{bibthbk}

\begin{thebibliography}{CGK06}

\bibitem[ABR96]{ande:brot:96}
L.~Andersen and R.~Brotherton-Ratcliffe.
\newblock Exact exotics.
\newblock {\em Risk}, 9:85--89, 1996.

\bibitem[AGP95]{asmu:glyn:pitm:95}
S.~Asmussen, P.~Glynn, and J.~Pitman.
\newblock Discretization error in simulation of one-dimensional reflecting
  {B}rownian motion.
\newblock {\em Ann. Appl. Probab.}, 5(4):875--896, 1995.

\bibitem[Als94]{alsm:94}
G.~Alsmeyer.
\newblock On the {M}arkov renewal theorem.
\newblock {\em Stoch. Proc. Appl.}, 50(1):37--56, 1994.

\bibitem[Avi07]{avik:07}
R.~Avikainen.
\newblock Convergence rates for approximations of functionals of {SDE}s.
\newblock {\em Submitted, available on http://arxiv.org/abs/0712.3635}, 2007.

\bibitem[Bal95]{bald:95}
P.~Baldi.
\newblock Exact asymptotics for the probability of exit from a domain and
  applications to simulation.
\newblock {\em Ann. Prob.}, 23(4):1644--1670, 1995.

\bibitem[BGK97]{broa:glas:kou:97}
M.~Broadie, P.~Glasserman, and S.~Kou.
\newblock A continuity correction for discrete barrier options.
\newblock {\em Mathematical Finance}, 7:325--349, 1997.

\bibitem[BL94]{boyl:lau:94}
P.P. Boyle and S.H. Lau.
\newblock Bumping up against the barrier with the binomial method.
\newblock {\em Jour. of Derivat.}, 1:6--14, 1994.

\bibitem[CGK06]{cost:elka:gobe:06}
C.~Costantini, E.~Gobet, and N.~El Karoui.
\newblock Boundary sensitivities for diffusion processes in time dependent
  domains.
\newblock {\em Appl. {M}ath. {O}ptim.}, 54--2:159--187, 2006.

\bibitem[DK06]{done:kypr:06}
R.~A. Doney and A.~E. Kyprianou.
\newblock Overshoots and undershoots of {L}\'evy processes.
\newblock {\em Ann. Appl. Probab.}, 16--1:91--106, 2006.

\bibitem[DL06]{deac:leja:06}
M.~Deaconu and A.~Lejay.
\newblock A random walk on rectangles algorithm.
\newblock {\em Methodology And Computing In Applied Probability},
  8(1):135--151, 2006.

\bibitem[FL01]{fuh:lai:01}
C.D. Fuh and T.L. Lai.
\newblock Asymptotic expansions in multidimensional {M}arkov renewal theory and
  first passage times for {M}arkov random walks.
\newblock {\em Adv. in Appl. Probab.}, 33(3):652--673, 2001.

\bibitem[Fre85]{frei:85}
M.~Freidlin.
\newblock {\em Functional integration and {P}artial differential equations}.
\newblock Annals of {M}athematics studies, {P}rinceton {U}niversity Press,
  1985.

\bibitem[Fri64]{frie:64}
A.~Friedman.
\newblock {\em Partial differential equations of parabolic type}.
\newblock Prentice-Hall, 1964.

\bibitem[GM04]{gobe:meno:spa:04}
E.~Gobet and S.~Menozzi.
\newblock Exact approximation rate of killed hypoelliptic diffusions using the
  discrete {E}uler scheme.
\newblock {\em Stoch. Proc. and Appl.}, 112:210--223, 2004.

\bibitem[GM05]{gobe:mair:05}
E.~Gobet and S.~Maire.
\newblock Sequential {C}ontrol {V}ariates for {F}unctionals of {M}arkov
  {P}rocesses.
\newblock {\em Siam. Journal of Num. Analysis}, 43-3:1256--1275, 2005.

\bibitem[GM07]{gobe:meno:semi:07}
E.~Gobet and S.~Menozzi.
\newblock Discrete sampling of functionals of {I}t\^o processes.
\newblock {\em S\'eminaire de Probabilit\'es}, XL:355--375, 2007.

\bibitem[Gob00]{gobe:00}
E.~Gobet.
\newblock Euler schemes for the weak approximation of killed diffusion.
\newblock {\em Stoch. Proc. Appl.}, 87:167--197, 2000.

\bibitem[Gob09]{gobe:09}
E.~Gobet.
\newblock {\em Handbook of Numerical Analysis, Vol. XV, Special Volume:
  Mathematical Modeling and Numerical Methods in Finance}, chapter Advanced
  Monte Carlo methods for barrier and related exotic options, pages 497--528.
\newblock Elsevier, Netherlands: North-Holland, 2009.

\bibitem[GT98]{gilb:trud:98}
D.~Gilbarg and N.S. Trudinger.
\newblock {\em Elliptic partial differential equations of second order}.
\newblock Springer Verlag, 1998.

\bibitem[Lie96]{lieb:96}
G.M. Lieberman.
\newblock {\em Second Order parabolic differential equations, 1st edn}.
\newblock World Scientfic, River Edge, NJ, 1996.

\bibitem[Men06]{meno:06}
S.~Menozzi.
\newblock Improved simulation for the killed {B}rownian motion in a cone.
\newblock {\em Siam Jour. Num. Anal.}, 44-6:2610--2632, 2006.

\bibitem[Mil97]{mils:97}
G.~N. Milstein.
\newblock Weak approximation of a diffusion process in a bounded domain.
\newblock {\em Stoch. Stoch. Reports}, 64:211--233, 1997.

\bibitem[MT99]{mils:tret:99}
G.~N. Milstein and M.~V. Tretyakov.
\newblock Simulation of a space-time bounded diffusion.
\newblock {\em Ann. Appl. Prob.}, 9--3:732--779, 1999.

\bibitem[RY99]{revu:yor:99}
D.~Revuz and M.~Yor.
\newblock {\em Continuous martingales and Brownian motion. 3rd ed.}
\newblock Grundlehren der Mathematischen Wissenschaften. 293. Berlin: Springer,
  1999.

\bibitem[Sie79]{sieg:79}
D.~Siegmund.
\newblock Corrected diffusion approximations in certain random walk problems.
\newblock {\em Adv. in Appl. Probab.}, 11(4):701--719, 1979.

\bibitem[Sie85]{sieg:85}
D.~Siegmund.
\newblock {\em Sequential Analysis}.
\newblock Springer, 1985.

\bibitem[Sim80]{simo:80}
J~Simon.
\newblock Differentiation with respect to the domain in boundary value
  problems.
\newblock {\em Numer. Funct. Anal. Optim.}, 7--8:649--687, 1980.

\bibitem[SY82]{sieg:yuh:82}
D.~Siegmund and Y.S. Yuh.
\newblock Brownian approximations for first passage probabilities.
\newblock {\em Z. Wahrsch. verw. Gebiete}, 59:239--248, 1982.

\bibitem[TT90]{tala:tuba:90}
D.~Talay and L.~Tubaro.
\newblock Expansion of the global error for numerical schemes solving
  sto\-chastic differential equations.
\newblock {\em Stoch. Anal. and App.}, 8-4:94--120, 1990.

\bibitem[Woo82]{wood:82}
M.~Woodroofe.
\newblock {\em Nonlinear renewal theory in sequential analysis}.
\newblock {SIAM}, 1982.

\bibitem[Zha88]{zhan:88}
C.H. Zhang.
\newblock A nonlinear renewal theory.
\newblock {\em Ann. Prob.}, 16-2:793--824, 1988.

\end{thebibliography}
\end{document}